# SPATIALLY EXPLICIT NON-MENDELIAN DIPLOID MODEL

By N. Lanchier and C. Neuhauser

*University of Minnesota*

We introduce a spatially explicit model for the competition between type $a$ and type $b$ alleles. Each vertex of the $d$-dimensional integer lattice is occupied by a diploid individual, which is in one of three possible states or genotypes: $aa$, $ab$ or $bb$. We are interested in the long-term behavior of the gene frequencies when Mendel's law of segregation does not hold. This results in a voter type model depending on four parameters; each of these parameters measures the strength of competition between genes during meiosis. We prove that with or without a spatial structure, type $a$ and type $b$ alleles coexist at equilibrium when homozygotes are poor competitors. The inclusion of a spatial structure, however, reduces the parameter region where coexistence occurs.

**1. Introduction.** Diploid organisms carry two sets of chromosomes, which are segregated during meiosis to produce haploid gametes. Mendel's law of segregation states that pairs of homologous chromosomes segregate randomly during meiosis and are distributed to gametes. This implies that a gamete of an individual who is heterozygous at a locus will receive one or the other allele with equal probability. This implies, for instance, that 25% of offspring of two individuals that are both heterozygous for the same pair of alleles will be homozygous for one allele, 25% will be homozygous for the other allele and 50% will be heterozygous.

This law does not always hold. Some organisms have chromosomal elements that show a higher than expected rate of transmission. A widely-studied example of meiotic drive can be found in the house mouse (*Mus musculus*). The house mouse may have a variant form of chromosome 17, called $t$-haplotype, which is transmitted to more than 90% of its offspring. Such systems that show a higher than 50% representation in their offspring are called meiotic drive systems.

Simple mathematical models readily show that meiotic drive by itself results in fixation. Yet natural populations show a low frequency of the $t$-haplotype [10–15%, see Durand et al. (1997) and references therein]. Two factors are contributing to the lower than expected frequency: male offspring who are homozygous for the $t$-haplotype are sterile (which by itself precludes fixation) and most $t$-haplotypes carry additional recessive lethal mutants. As Bruck (1957) showed, however, these two factors do not sufficiently explain the observed low frequencies. Meiotic drive has been observed in a number of other systems [see articles in Lyttle and Perkins (1991)], including the fruit fly Drosophila, mosquitoes, the fungus Neurospora, humans and maize. In order to detect meiotic drive in natural populations, a polymorphism must be present.

Lewontin and Dunn (1960) used computer simulations to suggest that random genetic drift due to small population size lowers Bruck's predicted frequency of the $t$-haplotype further. Other simulations with varying assumptions were conducted over the years [Levin, Petras and Rasmussen (1969); Petras and Topping (1983); Nunney and Baker (1993)], which are reviewed in Durand et al. (1997). Durand et al. (1997) concluded that "[N]one of the work surveyed here has fully explained the existence of the low frequency $t$-polymorphism seen in nature." They then posit a number of hypotheses to explain the low observed frequency, including low levels of migration in a spatially structured population, selection against the $t$-haplotype and lower transmission ratio distortion in wild versus laboratory populations. To test these hypotheses, they built a spatially structured simulation model. They found that "migration rate is the single most important factor in determining the equilibrium frequency of $t$-haplotypes." However, the model only yielded observed low frequencies for a very narrow range of parameters, and therefore essentially failed to provide an explanation for the low frequency of $t$-haplotypes in the wild. Their extensive simulations of this spatially structured model led this group to suggest that "the basic concept of a stable, low frequency $t$-polymorphism may be wrong." Instead, there might be two stable states, one in which the $t$-haplotype is absent and another where it is at high frequency (as in unstructured populations). Environmental change may be responsible for populations switching between the two states. They concluded with "[I]f these [transitions] occur frequently, many populations moving towards a new equilibrium (either $t$-haplotype extinction or panmictic frequencies) will be observed. The allele frequency measured in such populations would be low but transient."

Spatial simulations are difficult to interpret and we would like to caution against drawing conclusions that are not supported by rigorous mathematical analysis. Rigorous results are typically restricted to specific models, but generalities have emerged that allow us to speculate with some confidence on the behavior of models. While simulations might indicate the absence of a stable polymorphism for realistic parameter values and instead a presence of multiple stable states, we propose that a spatially explicit model on meiotic drive will either have a locally stable equilibrium, that is, a nontrivial stationary distribution leading to a stable polymorphism, or fixation will occur, where the fixation states correspond to the multiple stable states of the nonspatial model. However, we suggest that the model will not alternate among the fixation states (except in degenerate cases). This assertion is based on the observation that if a nonspatial model exhibits multiple, locally stable equilibria, it can often be shown that the corresponding spatially explicit, stochastic model will eventually converge to only one of these states and that all other fixation states are metastable [see, for instance, Neuhauser (1994)]. In addition, one- and two-dimensional spatial stochastic models often exhibit clustering. These clusters are not fixed in space but move around, and thus may give the impression of local transience. The model we now introduce appears to exhibit these commonly observed behaviors as well.

Our stochastic model is a continuous time Markov process in which the state at time $t$ is a function $\eta_t : \mathbb{Z}^d \longrightarrow \{aa, ab, bb\}$. Each site of the $d$-dimensional integer lattice is occupied by an individual with $\eta_t(x)$ indicating the genotype ($aa$, $ab$ or $bb$) of the individual present at site $x$ at time $t$. Each of the two genes at site $x$ gives birth at varying rates to a gene of its own type that then replaces one of the $4d$ genes picked at random among the $2d$ nearest neighbors. More precisely, the evolution at site $x$ is described by the transition rates

$$aa \to ab \quad \text{at rate } 2\phi_{bb} \sum_{\|x-z\|=1} \mathbb{1}\{\eta(z) = bb\} + \phi_{ba} \sum_{\|x-z\|=1} \mathbb{1}\{\eta(z) = ab\},$$

$$ab \to bb \quad \text{at rate } \phi_{bb} \sum_{\|x-z\|=1} \mathbb{1}\{\eta(z) = bb\} + \frac{\phi_{ba}}{2} \sum_{\|x-z\|=1} \mathbb{1}\{\eta(z) = ab\},$$

$$bb \to ab \quad \text{at rate } 2\phi_{aa} \sum_{\|x-z\|=1} \mathbb{1}\{\eta(z) = aa\} + \phi_{ab} \sum_{\|x-z\|=1} \mathbb{1}\{\eta(z) = ab\},$$

$$ab \to aa \quad \text{at rate } \phi_{aa} \sum_{\|x-z\|=1} \mathbb{1}\{\eta(z) = aa\} + \frac{\phi_{ab}}{2} \sum_{\|x-z\|=1} \mathbb{1}\{\eta(z) = ab\},$$

where $\|x - z\| = \sup_{i=1,2,\ldots,d} |x_i - z_i|$. For $i, j \in \{a, b\}$, the parameter $\phi_{ij}$ is the birth rate of alleles of type $i$ present in sites occupied by individuals of genotype $ij$. To understand the first two rates, note that an individual of genotype $bb$ produces genes of type $b$ at rate $2\phi_{bb}$ while an individual of genotype $ab$ produces genes of type $b$ at rate $\phi_{ba}$. If the gene (of type $b$) is sent to an individual of genotype $aa$, then it will replace one of the type $a$ genes which always results in a transition $aa \to ab$. But if the gene (of type $b$) is sent to an individual of genotype $ab$, then it will replace one of the genes chosen at random which results in a transition $ab \to bb$ with probability $1/2$. The last two rates can be understood similarly by exchanging the roles of $a$ and $b$.

*The mean-field model.* Before investigating the consequences of the inclusion of a spatial structure in the form of local interactions, we will look at the mean-field model [see Durrett and Levin (1994)]. The mean-field model can be obtained from the spatial model by assuming that individuals are located on the set of sites of a complete graph with $N$ vertices (i.e., each site has all other sites as its neighbors which implies that there is no spatial structure), by dividing the parameters $\phi_{ij}$ by the neighborhood size $N$ and by taking the limit as $N$ tends to infinity. This results in a nonspatial deterministic model in which the densities of genotypes evolve according to a system of differential equations. More precisely, by letting $u_{aa}$, $u_{ab}$ and $u_{bb}$ denote the densities of individuals of genotype $aa$, $ab$ and $bb$, respectively, the mean-field model is described by the following coupled system of ordinary differential equations:

$$\frac{du_{aa}}{dt} = (2\phi_{aa}u_{aa} + \phi_{ab}u_{ab})u_{ab}/2 - (2\phi_{bb}u_{bb} + \phi_{ba}u_{ab})u_{aa}, \tag{1}$$

$$\frac{du_{ab}}{dt} = (2\phi_{aa}u_{aa} + \phi_{ab}u_{ab})(u_{bb} - u_{ab}/2)$$
$$\tag{2}$$
$$\qquad\quad + (2\phi_{bb}u_{bb} + \phi_{ba}u_{ab})(u_{aa} - u_{ab}/2),$$

$$\frac{du_{bb}}{dt} = (2\phi_{bb}u_{bb} + \phi_{ba}u_{ab})u_{ab}/2 - (2\phi_{aa}u_{aa} + \phi_{ab}u_{ab})u_{bb}. \tag{3}$$

Let $\psi = (\phi_{aa} - \phi_{ba})(\phi_{bb} - \phi_{ab})$. When $\psi \leq 0$, the system has exactly two fixed points, the two trivial equilibria $u_{aa} = 1$ and $u_{bb} = 1$. By assuming that $\psi > 0$, we find an additional (nontrivial) equilibrium given by

$$\phi^2 u_{aa} = (\phi_{bb} - \phi_{ab})^2,$$
$$\phi^2 u_{bb} = (\phi_{aa} - \phi_{ba})^2,$$
$$\phi^2 u_{ab} = 2(\phi_{aa} - \phi_{ba})(\phi_{bb} - \phi_{ab}),$$

where $\phi = (\phi_{aa} - \phi_{ba}) + (\phi_{bb} - \phi_{ab})$. In Section 2 we will prove the following results:

1. The trivial equilibrium $u_{aa} = 1$ is locally stable if and only if $\phi_{aa} > \phi_{ba}$.
2. By symmetry, the equilibrium $u_{bb} = 1$ is locally stable if and only if $\phi_{bb} > \phi_{ab}$.
3. The interior fixed point ($\psi > 0$) is locally stable if and only if $\phi_{aa} < \phi_{ba}$ and $\phi_{bb} < \phi_{ab}$.

See Figure 1 for a summary of these results.

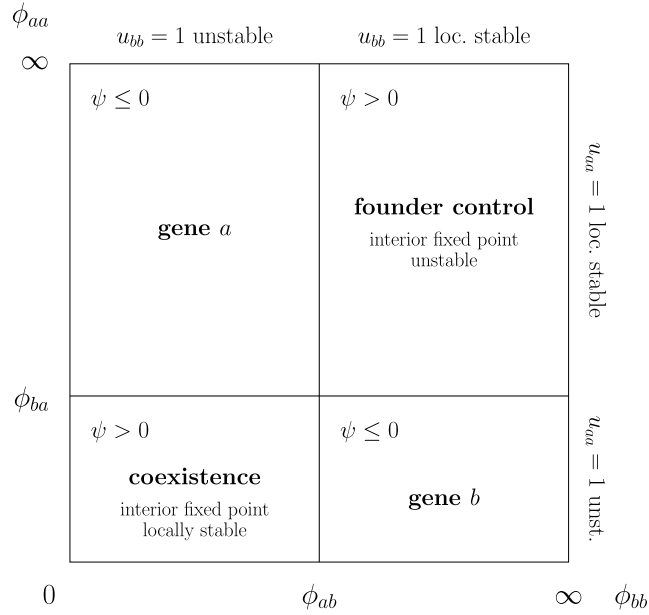

FIG. 1. *Phase diagram of the mean-field model. The vertical and horizontal axes represent the possible values of the rates $\phi_{aa}$ and $\phi_{bb}$, respectively, ranging from zero to infinity. In our picture, the rates $\phi_{ab}$ and $\phi_{ba}$ are constants whose fixed values are indicated at some arbitrary points.*

*The spatially explicit model.* To understand the role of spatial structure caused by local interactions, we now return to the spatially explicit model and provide comparisons with the nonspatial model. We refer the reader to Figures 1 and 2 that show the phase diagrams of the mean-field model and the interacting particle system, respectively. In Figure 1, the region labelled *gene a* indicates the set of parameters for which the density $u_{aa}$ converges to 1, the region *coexistence* the set of parameters for which there is a locally stable interior fixed point, and the region *founder control* the set of parameters for which the interior fixed point is unstable. In Figure 2, the region labelled *gene a* indicates the set of parameters for which the process starting with infinitely many individuals of type $aa$ converges to the "all $aa$" configuration and the region *coexistence* the set of parameters for which there is a stationary distribution with a positive density of both genes.

In the case when $\phi_{aa}$ and $\phi_{bb}$ are small, the particle system has a stationary distribution with a positive density of genes of type $a$ and $b$, which is symptomatic of the existence of a locally stable interior fixed point for the mean-field model when $\phi_{aa} < \phi_{ba}$ and $\phi_{bb} < \phi_{ab}$. This indicates that coexistence occurs for both models.

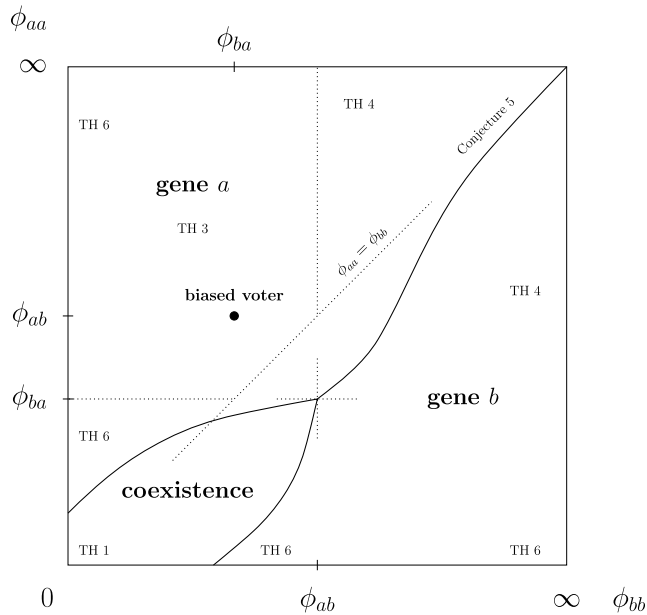

Fig. 2. *Spatially explicit model.*

On the other hand, in the case when $\phi_{ab}$ and $\phi_{ba}$ are small, coexistence is not possible. This, however, translates differently for the spatial and non-spatial models. In the spatial model, we believe that, for almost all the parameters, there is fixation of one of the genes. Provided that the process starts with infinitely many genes of both types, the "winner" depends on the parameters of the model (see Figure 2). Simulations suggest that, along the transition curve, clustering occurs in any dimension (see the right-hand side of Figure 3 for a picture in $d = 2$), and that the three regimes meet at point $\phi_{aa} = \phi_{ba}$ and $\phi_{bb} = \phi_{ab}$ as indicated in Figure 2. In contrast, the mean-field model exhibits founder control when $\phi_{aa} > \phi_{ba}$ and $\phi_{bb} > \phi_{ab}$, that is, the interior fixed point is unstable and the "winner" depends on both the parameters and the initial densities (see the right-hand side of Figure 4).

Numerical simulations also indicate that the set of parameters for which coexistence may occur is sensibly smaller for the spatial model. We will prove rigorously the existence of a set of parameters for which the interior fixed point is locally stable for the mean-field model whereas coexistence is not possible for the 1-dimensional stochastic model. This reduction of the coexistence region of the particle system when compared to the mean-field model has already been proved theoretically for a spatially explicit version of the Lotka–Volterra model [see Neuhauser and Pacala (1999)]. We now state our results in detail starting with the coexistence result.

THEOREM 1 (Coexistence). *Assume that $\phi_{ab} = \phi_{ba} > 0$. If $\phi_{aa}$ and $\phi_{bb}$ are sufficiently small then there exists a translation invariant stationary distribution $\nu$ such that $\nu(\eta(x) = ab) > 0$. The same holds when $\phi_{ab} < 2\phi_{ba}$ and $\phi_{ba} < 2\phi_{ab}$ in the 1-dimensional case.*

See the left-hand side of Figure 3 for a picture. In the case when $\phi_{ab} = \phi_{ba} = 1$ and $\phi_{aa} = \phi_{bb} = 0$, the process $\eta_t$ reduces to an annihilating branching process (ABP) in which 0 = homozygote of either type and 1 = heterozygote, that is, the state of each site of the lattice flips from 0 to 1 and from 1 to 0 at rate equal to the number of neighbors in state 1. The main ingredient to prove Theorem 1 is a rescaling argument that compares the ABP viewed on suitable length and time scales with a 1-dependent oriented percolation process on the lattice

$$\mathcal{G} = \{(z, n) \in \mathbb{Z}^2 : z + n \text{ is even and } n \geq 0\}.$$

More precisely, by defining occupied sites in the lattice $\mathcal{G}$ as the set of sites $(z, n) \in \mathcal{G}$ such that there is at least one particle (state 1) in the spatial cube

$$Nze_1 + [-N/4, N/4]^d = (Nz, 0, \ldots, 0) + [-N/4, N/4]^d$$

at time $nT$, we will prove the following theorem.

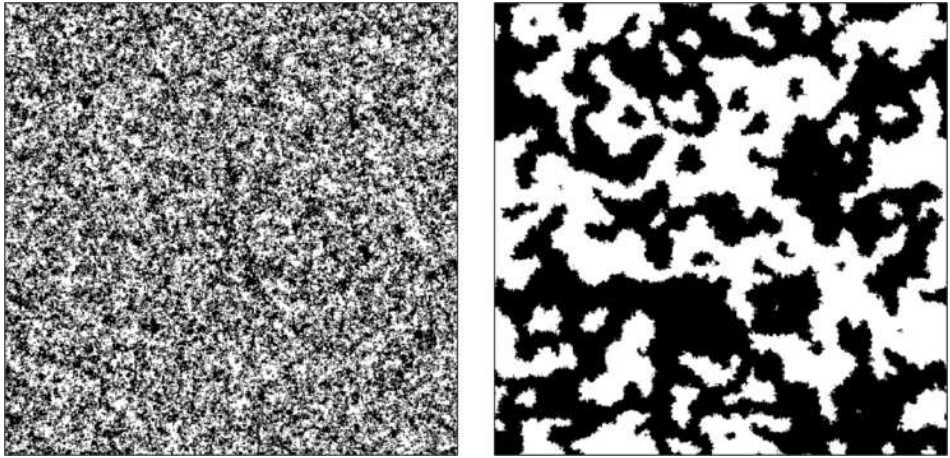

FIG. 3. *Spatially explicit model with nearest neighbor interactions at time 50 on the $400 \times 400$ square with periodic boundary conditions. White (respectively, black) sites refer to homozygotes of type a (respectively, b), and grey sites to heterozygotes. Left: $\phi_{aa} = \phi_{bb} = 4$ and $\phi_{ab} = \phi_{ba} = 5$. Right: $\phi_{aa} = \phi_{bb} = 5$ and $\phi_{ab} = \phi_{ba} = 1$. In both pictures, the process starts with a Bernoulli product measure, half of the genes being of type a. The picture on the left shows the process at equilibrium but not the picture on the right in which clusters grow in time indicating that the only equilibria are the "all aa" and the "all bb" configurations.*

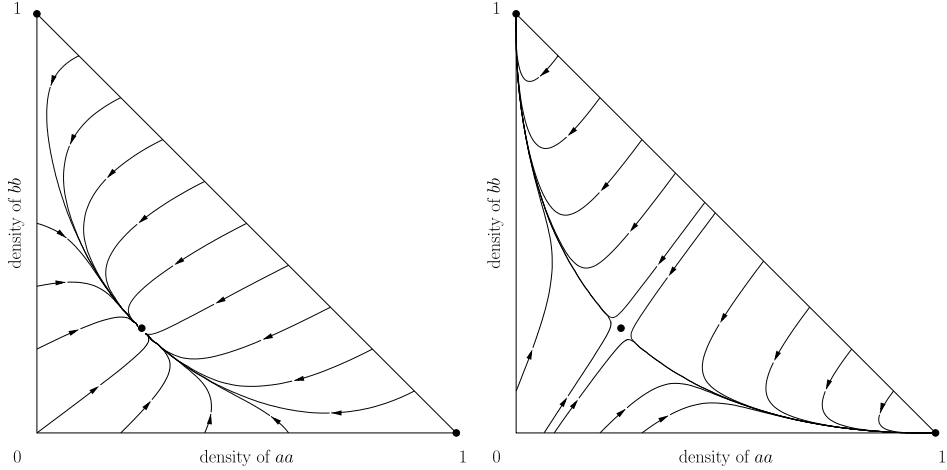

Fig. 4. *Solution curves of the mean-field model. The parameters are given by $\phi_{aa} = 1$, $\phi_{ba} = 3$, $\phi_{bb} = 2$ and $\phi_{ab} = 4$ (left), and $\phi_{aa} = 4$, $\phi_{ba} = 2$, $\phi_{bb} = 3$ and $\phi_{ab} = 1$ (right).*

THEOREM 2 (Annihilating branching process). *For any $\varepsilon > 0$, the parameters $N$ and $T$ can be chosen in such a way that the set of occupied sites dominates the set of wet sites in a 1-dependent oriented site percolation process on $\mathcal{G}$ in which sites are open with probability $1 - \varepsilon$.*

This implies coexistence of both genes (survival of the heterozygotes) when $\phi_{aa} = \phi_{bb} = 0$. To conclude the proof of Theorem 1, we will apply a standard perturbation argument to extend the result to a region where $\phi_{aa}$ and $\phi_{bb}$ are positive but small.

The last three theorems provide conditions under which there is fixation of the gene of type $a$ and so extinction of the gene of type $b$. Here, "fixation of the gene of type $a$" means that the process converges to the "all $aa$" configuration. The reader will note that, by symmetry of the evolution rules of the process, the conditions under which there is fixation of the gene of type $b$ can be deduced by reversing the roles of genes of type $a$ and $b$. The following results hold for the process starting with infinitely many genes of type $a$ or infinitely many individuals of genotype $aa$. We can weaken this condition and assume instead that the process starts with finitely many genes of type $a$ or finitely many individuals of genotype $aa$, and conclude that fixation of the gene of type $a$ occurs with positive probability.

THEOREM 3 (Fixation). *Assume that $\min(\phi_{aa}, \phi_{ab}) > \max(\phi_{ba}, \phi_{bb})$ and that infinitely many genes of type $a$ are present at time 0. Then there is fixation of the gene of type $a$.*

THEOREM 4 (Fixation). *Assume that $\phi_{ab}$ is fixed,*

$$\phi_{aa} > 6d^2 \left(1 + \sqrt{1 + \frac{2}{d}}\right) \phi_{bb} \tag{4}$$

*and that infinitely many individuals of genotype aa are present at time 0. Then there is fixation of the gene of type a provided $\phi_{ba}$ is sufficiently small.*

By assuming that $\phi_{aa} = \phi_{ab}$ and $\phi_{ba} = \phi_{bb}$, the birth rate of a gene (of either type) at site $x$ does not depend upon the type of the second gene present at site $x$. This makes our process a biased voter model and allows us to define naturally a dual process. In all the other cases, duality is not tractable. The proofs of Theorems 3 and 4 are made difficult by this lack of a dual process. To establish Theorem 3, we will rely on a coupling argument to compare the process with a biased voter model in which genes of type $a$ are more competitive than genes of type $b$. The proof of Theorem 4 is more difficult. The idea is to prove that, under condition (4) and for the process starting with all sites in $F_N = [-N, N]^d$ in state $aa$, the probability that genes of type $b$ reach site 0 decreases exponentially fast with the length side $N$. This will be estimated by studying the evolution of the process along all the paths starting outside $F_N$ and ending at site 0. Then we will use the individual of genotype $aa$ at site 0 as a source to produce a large cube void of genes of type $b$ and will conclude by comparing the process viewed on suitable length and time scales with oriented percolation. A time rescaling in Theorem 4 suggests that $\phi_{ba}$ being fixed if the parameters $\phi_{aa}$ and $\phi_{bb}$ are sufficiently large, and if condition (4) holds, then there is fixation of the gene of type $a$. Unfortunately, this cannot be deduced from our result since $\phi_{ba}$ in the statement of Theorem 4 has to be smaller than some constant that (a priori) depends upon the values of the parameters $\phi_{aa}$ and $\phi_{bb}$. However, numerical simulations suggest the following.

CONJECTURE 5 (Fixation). *Assume that $\phi_{ab}$ and $\phi_{ba}$ are fixed and that infinitely many genes of type $a$ are present at time 0. Also, assume that the ratio $\phi_{aa}/\phi_{bb} > 1$ is fixed. Then there is fixation of the gene of type $a$ provided $\phi_{aa}$ is sufficiently large.*

Numerical simulations also indicate that when $\phi_{aa} = \phi_{bb} > \phi_{ab} = \phi_{ba}$ (founder control for the mean-field model), clustering occurs in any dimension for the spatial model. The right-hand side of Figure 3 gives an illustration in $d = 2$.

Theorems 1 and 3 show similarities with the mean-field model. Theorem 4 implies the existence of a parameter region for which fixation occurs for the spatial model (the "winner" depends on the choice of the parameters) while the mean-field model exhibits founder control. In both models, coexistence cannot occur for these parameters. Similarly, Conjecture 5 and Theorem 6, formula (5), below which confirms our conjecture in $d=1$ extend the parameter region for which one of the genes goes extinct eventually. On the contrary, Theorem 6, formula (6), shows a major difference between spatial and nonspatial models, namely the set of parameters for which coexistence occurs is smaller for the spatial model (we believe that this holds in any dimension, though our proof relies heavily on geometric arguments true in $d=1$ only). This might provide a theoretical explanation for the low frequency of the $t$-haplotype in the case of the house mouse. To model the evolution of the house mouse, we first choose $\phi_{ba} < \phi_{ab}$ where $a$ represents the $t$-haplotype to take into account the fact that the $t$-haplotype is transmitted to a majority of the offspring. Since male offspring who are homozygous for the $t$-haplotype are sterile, we also fix $\phi_{aa} = 0$, though our model does not distinguish between male and female. Finally, we let the last parameter be $\phi_{bb} = (1/2)(\phi_{ab} + \phi_{ba})$, assuming that $ab$ individuals are overall equally fit as $bb$ individuals. This set of parameters corresponds to points of the phase diagram which are well inside the coexistence phase for the nonspatial model (see Figure 1) but to points which are much closer to the boundary of the coexistence phase where the frequency of type $a$ should be expected to be lower for the spatial model (see Figure 2). In particular, the spatial component might be responsible for the low frequency of the $t$-haplotype in the example of the house mouse.

THEOREM 6 (Fixation). *Assume that $d=1$ and that infinitely many individuals of genotype $aa$ are present at time 0. Then whenever one the following two conditions holds:*

$$(5) \qquad \phi_{aa} > \phi_{bb} + \phi_{ba} + \sqrt{\phi_{bb}\phi_{ba}},$$

$$(6) \qquad \phi_{aa} > \max\left(\frac{2\phi_{ba}}{5}, \phi_{ba} - \frac{\phi_{ab}}{6}\right) \quad \text{and} \quad \phi_{bb} > 0 \text{ small,}$$

*there is fixation of the gene of type $a$.*

The rest of this paper is devoted to proofs. In Section 2, we investigate in detail the existence and the stability of the fixed points of the mean-field equations (1)–(3). In Sections 3–6, we will prove the results regarding the spatially explicit model.

**2. The mean-field model.** This section is devoted to the study of the mean-field model introduced in Section 1. First of all, by setting the right-hand sides of (1) and (3) equal to 0, we obtain

$$2u_{aa}(2\phi_{bb}u_{bb} + \phi_{ba}u_{ab}) = u_{ab}(2\phi_{aa}u_{aa} + \phi_{ab}u_{ab}),$$
$$2u_{bb}(2\phi_{aa}u_{aa} + \phi_{ab}u_{ab}) = u_{ab}(2\phi_{bb}u_{bb} + \phi_{ba}u_{ab}),$$

so that $4u_{aa}u_{bb} = u_{ab}^2$. Note that the previous equation holds when $2\phi_{aa}u_{aa} + \phi_{ab}u_{ab} = 0$ as well since in this case $u_{aa} = u_{ab} = 0$. Substituting this into (1) and setting the right-hand side equal to 0, leads to

$$\frac{du_{aa}}{dt} = (\phi_{aa} - \phi_{ba})u_{aa}u_{ab} + (1/2)\phi_{ab}u_{ab}^2 - 2\phi_{bb}u_{aa}u_{bb}$$

$$= (\phi_{aa} - \phi_{ba})u_{aa}u_{ab} - 2(\phi_{bb} - \phi_{ab})u_{aa}\ u_{bb} = 0,$$

from which it follows that

$$u_{aa} = 0 \quad \text{or} \quad (\phi_{aa} - \phi_{ba})u_{ab} = 2(\phi_{bb} - \phi_{ab})u_{bb}$$

at equilibrium. Similarly, using (3), we find

$$u_{bb} = 0 \quad \text{or} \quad (\phi_{bb} - \phi_{ab})u_{ab} = 2(\phi_{aa} - \phi_{ba})u_{aa}$$

at equilibrium. The condition $u_{aa} = 0$ gives rise to the trivial equilibrium $u_{bb} = 1$, and the condition $u_{bb} = 0$ to the trivial equilibrium $u_{aa} = 1$. If $u_{aa}u_{bb} \neq 0$, then

$$2(\phi_{aa} - \phi_{ba})^2 u_{aa} = (\phi_{aa} - \phi_{ba})(\phi_{bb} - \phi_{ab})u_{ab} = 2(\phi_{bb} - \phi_{ab})^2 u_{bb}.$$

Setting $\phi = (\phi_{aa} - \phi_{ba}) + (\phi_{bb} - \phi_{ab})$, we get

$$\phi^2 u_{aa} = (\phi_{aa} - \phi_{ba})^2 u_{aa} + 2(\phi_{aa} - \phi_{ba})(\phi_{bb} - \phi_{ab})u_{aa} + (\phi_{bb} - \phi_{ab})^2 u_{aa}$$

$$= (\phi_{bb} - \phi_{ab})^2 u_{bb} + (\phi_{bb} - \phi_{ab})^2 u_{ab} + (\phi_{bb} - \phi_{ab})^2\ u_{aa}$$

$$= (\phi_{bb} - \phi_{ab})^2.$$

By using the symmetry of the model, we deduce that at equilibrium

$$\phi^2 u_{aa} = (\phi_{bb} - \phi_{ab})^2,$$
$$\phi^2 u_{bb} = (\phi_{aa} - \phi_{ba})^2,$$
$$\phi^2 u_{ab} = 2(\phi_{aa} - \phi_{ba})(\phi_{bb} - \phi_{ab}).$$

Therefore, the condition for the existence of a nontrivial equilibrium is given by

$$\psi := (\phi_{aa} - \phi_{ba}) \times (\phi_{bb} - \phi_{ab}) > 0.$$

The stability of the equilibria can be analyzed using standard linearization techniques. Namely, denoting by $\mathcal{J}_{eq}$ the Jacobian matrix at $eq$, the equilibrium is locally stable if and only if all the eigenvalues of $\mathcal{J}_{eq}$ associated with eigenvectors $v = (v_1, v_2, v_3)$ such that $v_1 + v_2 + v_3 = 0$ are negative. We first look at the stability of $u_{aa} = 1$. The Jacobian matrix is given by

$$\mathcal{J}_a = \begin{pmatrix} 0 & \phi_{aa} - \phi_{ba} & -2\phi_{bb} \\ 0 & -\phi_{aa} + \phi_{ba} & 2\phi_{aa} + 2\phi_{bb} \\ 0 & 0 & -2\phi_{aa} \end{pmatrix}.$$

The eigenvalues of the matrix $\mathcal{J}_a$ are 0, $-\phi_{aa}+\phi_{ba}$ and $-2\phi_{aa}$. The eigenspaces associated with the eigenvalues 0 and $-\phi_{aa}+\phi_{ba}$ are respectively, given by

$$\text{Span}((1,0,0)) \quad \text{and} \quad \text{Span}((1,-1,0))$$

so that the equilibrium $u_{aa}=1$ is locally stable if and only if $\phi_{aa}>\phi_{ba}$. By symmetry, the equilibrium $u_{bb}=1$ is locally stable if and only if $\phi_{bb}>\phi_{ab}$. Now, by letting

$$\psi_{aa} = (\phi_{aa}-\phi_{ba})^2,$$
$$\psi_{bb} = (\phi_{bb}-\phi_{ab})^2,$$
$$\psi_{ab} = (\phi_{aa}-\phi_{ba})(\phi_{bb}-\phi_{ab}),$$

we find that the Jacobian matrix at the interior fixed point is given by

$$\mathcal{J}_{ab} = \frac{1}{\phi^2} \begin{pmatrix} -2\psi_{aa}\phi_{ab} & \psi_{ab}(\phi_{bb}+\phi_{ab}) & -2\psi_{bb}\phi_{bb} \\ 2\psi_{aa}(\phi_{aa}+\phi_{ab}) & -\psi_{ab}(\phi_{aa}+\phi_{bb}+\phi_{ab}+\phi_{ba}) & 2\psi_{bb}(\phi_{bb}+\phi_{ba}) \\ -2\psi_{aa}\phi_{aa} & \psi_{ab}(\phi_{aa}+\phi_{ba}) & -2\psi_{bb}\phi_{ba} \end{pmatrix},$$

where $\phi = (\phi_{aa}-\phi_{ba})+(\phi_{bb}-\phi_{ab})$ as above. First of all, it is easy to check that 0 is an eigenvalue, and that it is associated to the eigenvector

$$((\phi_{bb}-\phi_{ab})^2, 2(\phi_{aa}-\phi_{ba})(\phi_{bb}-\phi_{ab}), (\phi_{aa}-\phi_{ba})^2)$$

whose coordinates add up to 0 only if $\psi \leq 0$. In particular, since the condition $\psi > 0$ is required for the existence of an interior fixed point, the eigenvalue 0 is irrelevant in the analysis of the stability. Now, we introduce the curve

$$\Gamma = \{(u_{aa}, u_{ab}, u_{bb}): u_{aa}=p^2, u_{ab}=2p(1-p),$$
$$u_{bb}=(1-p)^2 \text{ for some } p \in [0,1]\}.$$

This curve is clearly visible in Figure 4. A population of diploid individuals in which the genotype frequencies lie on $\Gamma$ is said to be in Hardy–Weinberg equilibrium. When Mendel's law of segregation holds, points in this curve correspond to the limiting genotype frequencies of a nonspatial population (mean-field model), with $p$ denoting the initial frequency of genes of type $a$. Note however that the inclusion of local interactions translates into spatial correlations that lead to an increase of the density of homozygous individuals, and a decrease of the density of heterozygous individuals. For the general mean-field model (when Mendel's law does not hold), points in $\Gamma$ are no longer equilibrium points, but the curve is globally invariant under the time evolution. Observing in addition that $\Gamma$ contains the interior fixed point (take $p=(\phi_{bb}-\phi_{ab})\phi^{-1}$), we find that

$$((\phi_{bb}-\phi_{ab}), (\phi_{aa}-\phi_{ba})-(\phi_{bb}-\phi_{ab}), -(\phi_{aa}-\phi_{ba}))$$

is an eigenvector, and that it is associated to the eigenvalue

$$\phi_1 := (\phi_{aa}-\phi_{ba})(\phi_{bb}-\phi_{ab})((\phi_{aa}-\phi_{ba})+(\phi_{bb}-\phi_{ab})).$$

This shows that, in the direction of the curve $\Gamma$, the interior fixed point is:

1. locally stable if $\phi_{aa} < \phi_{ba}$ and $\phi_{bb} < \phi_{ab}$,
2. unstable if $\phi_{aa} > \phi_{ba}$ and $\phi_{bb} > \phi_{ab}$.

This proves the instability when $\phi_{aa} > \phi_{ba}$ and $\phi_{bb} > \phi_{ab}$. See Figure 4 for a picture of the solution curves. Finally, using that 0 and $\phi_1$ are eigenvalues, it is straightforward to conclude that (as the third root of the characteristic polynomial) the third eigenvalue is given by

$$\phi_2 := -2((\phi_{aa} - \phi_{ba})(\phi_{bb} - \phi_{ab})(\phi_{aa} + \phi_{bb})$$
$$+ (\phi_{bb} - \phi_{ab})^2 \phi_{ba} + (\phi_{aa} - \phi_{ba})^2 \phi_{ab})$$

which is always negative whenever $\psi = (\phi_{aa} - \phi_{ba})(\phi_{bb} - \phi_{ab}) > 0$, the condition required for the existence of an interior fixed point. This proves the stability when $\phi_{aa} < \phi_{ba}$ and $\phi_{bb} < \phi_{ab}$.

**3. Proof of Theorem 1.** The aim of this section is to prove that both alleles may coexist provided $\phi_{aa}$ and $\phi_{bb}$ are sufficiently small. We will start by assuming that $\phi_{ab} = \phi_{ba} > 0$ and will explain at the end of this section how to deduce the result in the 1-dimensional case when $\phi_{ab} < 2\phi_{ba}$ and $\phi_{ba} < 2\phi_{ab}$.

First, we observe that when $\phi_{ab} = \phi_{ba} > 0$ and $\phi_{aa} = \phi_{bb} = 0$, the local dynamics of heterozygous sites can be defined regardless of the type of nearby homozygous sites. More precisely, by setting 0 = homozygote of either type and 1 = heterozygote, or, equivalently, by letting $\bar{\eta}_t$ denote the stochastic process defined for any $x \in \mathbb{Z}^d$ and any time $t \geq 0$ by

(7) $$\bar{\eta}_t(x) = \begin{cases} 1, & \text{if } \eta_t(x) = ab, \\ 0, & \text{otherwise,} \end{cases}$$

we obtain a new Markov process with state space $\{0,1\}^{\mathbb{Z}^d}$ whose state at each lattice site flips from 0 to 1 and from 1 to 0 at rate $\phi_{ab}$ times the number of neighbors in state 1. This process is known as the annihilating branching process (ABP) and has been investigated by Sudbury (1990) and Bramson, Ding and Durrett (1991). Their results show that starting with any measure, which puts no mass on the empty configuration, the ABP converges in distribution to the Bernoulli product measure with density 1/2. This implies that alleles of type $a$ and $b$ coexist when $\phi_{aa} = \phi_{bb} = 0$. To prove Theorem 1, we will employ a rescaling argument [see Bramson and Durrett (1988)] where we tile the space–time region $\mathbb{Z}^d \times \mathbb{Z}^+$ into large boxes, define a "good event" in each box and show that if a good event occurs in one such box, it will with probability close to 1 occur in neighboring boxes at a later time. Comparison with an oriented percolation process on $\mathbb{Z} \times \mathbb{Z}^+$ then completes the proof. The "good event" is the presence of heterozygous individuals in a given box if one thinks of the process $\eta_t$, or the presence of 1's in a given box if one thinks of the process $\bar\eta_t$. To show that heterozygosity can spread, we define a repositioning algorithm that shows that heterozygous individuals spread out linearly in the direction of any of the coordinate axes (Lemma 3.1). The reason for invoking the rescaling argument is that once the process is embedded in a block structure, it is straightforward to apply a perturbation argument and deduce that coexistence occurs as well when $\phi_{aa}$ and $\phi_{bb}$ are sufficiently small. The main ingredient to compare the ABP with a percolation process is that particles can be moved along a given straight line $\Delta$ while staying, with probability close to 1, within a reasonably short distance of $\Delta$ (see Lemma 3.1 below).

For $i \in \{1, 2, \ldots, d\}$, we let $e_i$ denote the $i$th unit vector, $\pi_i$ denote the orthogonal projection on the $i$th axis, and $H_i$ denote the hyperplane orthogonal to $e_i$, that is,

$$\pi_i(x) = \langle x, e_i \rangle = x_i \quad \text{and} \quad H_i = \mathrm{Span}(e_1, e_2, \ldots, e_{i-1}, e_{i+1}, \ldots, e_d),$$

where $\langle \cdot, \cdot \rangle$ is the usual inner product on $\mathbb{R}^d$. We also set

$$H_{i,k} = [-k^{0.7}, k^{0.7}]^d \cap H_i \quad \text{and} \quad E_{i,k} = \bigcup_{j=-k}^{k} (H_{i,k} + je_i),$$

where $k$ is a large integer whose value will vary (increase) from lemma to lemma. For $\omega \in \{-1, +1\}$ and $j = 1, 2, \ldots, k$, we introduce the hitting times

$$\tau_{\omega,i,j,k} = \inf\{t \geq 0 : \bar\eta_t(x) = 1 \text{ for some } x \in H_{i,k} + \omega j e_i\}$$

and $\tau_{\omega,i,k} = \sup_{j=1,2,\ldots,k} \tau_{\omega,i,j,k}$. Roughly speaking, the next lemma tells us that, with probability close to 1 when $k$ is large, a particle at site 0 can be moved of a distance $k$ along the $i$th axis while staying within a distance $k^{0.7}$ of this axis. This can be done in less than $Bk$ units of time for some constant $B$ that depends on the parameters of the process but not on $k$.

LEMMA 3.1. *Let $\bar\eta_0(0) = 1$ and $i \in \{1, 2, \ldots, d\}$. There exist $B < \infty$, $C_1 < \infty$ and $\alpha_1 > 0$ independent of $k$ such that*

$$P(\tau_{\omega,i,k} > Bk) \leq C_1 \exp(-\alpha_1 k^{0.1}) \quad \text{and}$$

$$P(\bar\eta_t(x) = 0 \text{ for all } x \in E_{i,k} \text{ and some } t \leq \tau_{\omega,i,k}) \leq C_1 \exp(-\alpha_1 k^{0.1})$$

*for $k$ sufficiently large.*

PROOF. Note that, by symmetry, it suffices to prove Lemma 3.1 for $\omega = +1$. The idea is to construct a jump process $X_t$ that keeps track of a particle moving towards the hypercube $H_{i,k} + k e_i$ while staying within a distance $k^{0.7}$ of the $i$th axis. The process $X_t$ starts at $X_0 = 0$ and the evolution rules are formulated according to whether the process moves in the direction of the vector $e_i$ or in the hyperplane $X_t + H_i$ as follows:

1. Whenever $\bar{\eta}_t(X_t + e_i) = 1$, the process $X_t$ jumps instantaneously to $X_t + e_i$. In particular, site $X_t + e_i$ is empty, that is, in state 0, at any time. If the particle at $X_t$ is killed by a particle located at $X_t - e_i$ (due to the birth of a gene originated from site $X_t - e_i$), then $X_t$ jumps instantaneously to $X_t - e_i$.
2. If the particle at site $X_t$ gives birth to a particle which is sent to $x \in X_t + H_i$, then $X_t$ jumps to site $x$. If the particle at site $X_t$ is killed by a particle located at $x \in X_t + H_i$ (due to the birth of a gene originated from site $x \in X_t + H_i$), then $X_t$ jumps to $x$.

Rules 1 and 2 are illustrated in Figure 5 (note that the definition of $X_t$ implies that at time $t$ site $X_t$ is always occupied). The picture shows a configuration of the process $\bar{\eta}_t$ together with the position of $X_t$. Open circles refer to empty sites (state 0) and full circles to occupied sites (state 1). If a particle gives birth through a $j$-arrow, then the process $X_t$ jumps to site $x_j$, $j = 1, 2, 3, 4$. Births through 1- and 2-arrows illustrate rule 1 (jumps in the direction of $e_i$), while births through 3- and 4-arrows illustrate rule 2 (jumps in the direction of $H_i$). Lemma 3.1 is a straightforward consequence of Lemmas 3.2 and 3.3 in which we investigate the projections of the process $X_t$ on the $i$th axis and the hyperplane $H_i$, respectively.

In the next lemma, we prove that rule 1 above implies that, with probability close to 1 for $k$ large, the projection $\pi_i(X_t)$ reaches $k$ in less than $Bk$ units of time and before reaching $-k$.

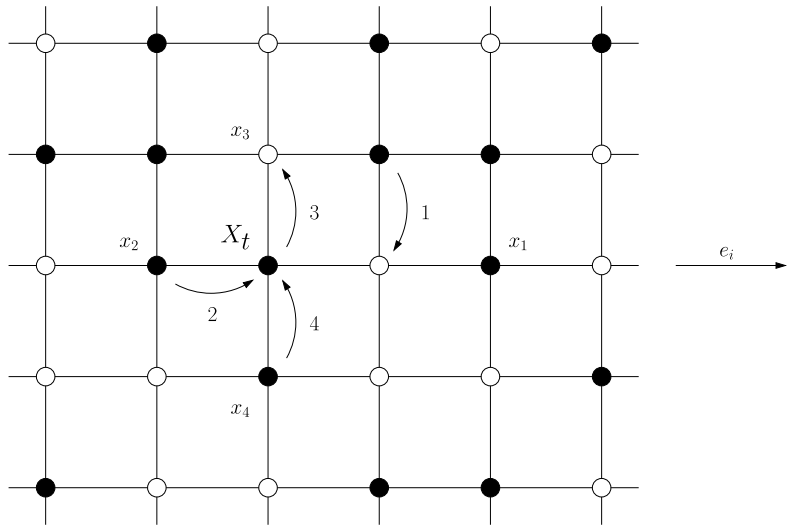

FIG. 5. *Evolution rules of the process $X_t$.*

LEMMA 3.2. *Let $\sigma_{\omega,i} = \inf\{t \geq 0 : \pi_i(X_t) \geq k\}$. Then there exists $B < \infty$ such that, for $k$ sufficiently large,*

$$P(\sigma_{\omega,i} > Bk) \leq C_2 \exp(-\alpha_2 k) \quad \text{and}$$

$$P(\pi_i(X_t) \leq -k \text{ for some } t \leq \sigma_{\omega,i}) \leq C_2 \exp(-\alpha_2 k)$$

*for suitable constants $C_2 < \infty$ and $\alpha_2 > 0$.*

PROOF. We observe that $\pi_i(X_t)$ increases by one whenever a particle is sent to site $X_t + e_i$ (type 1 arrow in Figure 5), which occurs at rate at least $\phi_{ab} = \phi_{ba}$ since $X_t$ is occupied. Rule 1 above also implies that the process $\pi_i(X_t)$ can decrease by one only if $X_t - e_i$ is occupied, the decrease resulting from the birth at $X_t$ of a gene that originated from site $X_t - e_i$ (type 2 arrow in Figure 5). This occurs (when $X_t - e_i$ is occupied) at rate $\phi_{ab}/2 + \phi_{ba}/2 = \phi_{ab}$. In other words, if we introduce the gap process

$$G_t = \inf\{j \geq 0 : \bar{\eta}_t(X_t - (j+1)e_i) = 1\},$$

then $\pi_i(X_t)$ jumps to the interval $[\pi_i(X_t) + 1, \infty)$ at rate at least $\phi_{ab}$ regardless of the value of the process $G_t$ and jumps to $\pi_i(X_t) - 1$ at rate $\phi_{ab}$ only when $G_t = 0$. In addition, a straightforward calculation shows that $G_t$ makes transitions as follows:

$$0 \to I_1 \quad \text{at rate at least } 2\phi_{ab},$$
$$I_1 \to 0 \quad \text{at rate at most } 4d\phi_{ab},$$

where $I_1 = \{1, 2, \ldots\}$ and where "$0 \to I_1$" (respectively, "$I_1 \to 0$") stands for some transition from 0 to $m$ (respectively, from $m$ to 0) for some $m \in I_1$. When site $X_t - 2e_i$ is empty, the first rate is the rate at which particles at $X_t - e_i$ and $X_t$ kill each other. When $X_t - 2e_i$ is occupied, it is the rate at which particles at $X_t - 2e_i$ and $X_t$ kill the particle at site $X_t - e_i$. The upper bound for the second rate is reached when all the neighbors of sites $X_t - e_i$ and $X_t + e_i$ are occupied. The evolution rules imply that

$$P(G_{t+h} = 0) \leq P(G_t = 0) \times [1 - 2\phi_{ab}h] + P(G_t \neq 0) \times 4d\phi_{ab}h + O(h^2).$$

By subtracting $P(G_t = 0)$ from both sides, dividing by $h$ and taking the limit as $h \to 0$,

$$\frac{d}{dt} P(G_t = 0) \leq -2\phi_{ab} \times P(G_t = 0) + 4d\phi_{ab} \times P(G_t \neq 0)$$
$$\leq 4d\phi_{ab} - (4d + 2)\phi_{ab} P(G_t = 0)$$

from which it follows that

$$P(G_t = 0) \leq \frac{2d}{2d+1} + \left[P(G_0 = 0) - \frac{2d}{2d+1}\right] \times \exp(-(4d+2)\phi_{ab} t).$$

Taking the lim sup as $t \to \infty$, we get

$$\limsup_{t \to \infty} P(G_t = 0) \leq \frac{2d}{2d+1}.$$

Finally, by using that $\pi_i(X_t)$ increases by one at rate at least $\phi_{ab}$, and decreases by one at rate at most $\phi_{ab}$ when $G_t = 0$ and at rate 0 when $G_t \neq 0$, we obtain

$$\liminf_{t \to \infty} \frac{\pi_i(X_t)}{t} \geq \phi_{ab} - \phi_{ab} \times \limsup_{t \to \infty} P(G_t = 0) \geq \frac{\phi_{ab}}{2d+1}. \tag{8}$$

In particular, by setting $B = 4d\phi_{ab}^{-1}$, we can conclude that

$$P(\sigma_{\omega,i} > Bk) \leq P(\pi_i(X_t) \leq k \text{ for all } t \leq Bk) \leq P(\pi_i(X_{Bk}) \leq k)$$
$$= P(\pi_i(X_{Bk}) \leq (\phi_{ab}/4d) \times Bk) \leq C_3 \exp(-\alpha_3 k)$$

for suitable $C_3 < \infty$ and $\alpha_3 > 0$ and $k$ sufficiently large. The limit in (8) also implies that

$$P(\pi_i(X_t) \leq -k \text{ for some } t \leq \sigma_{\omega,i}) \leq P(\pi_i(X_t) \leq -k \text{ for some } t \geq 0)$$
$$\leq C_4 \exp(-\alpha_4 k)$$

for suitable $C_4 < \infty$ and $\alpha_4 > 0$ and $k$ sufficiently large. The lemma follows. $\square$

The next lemma shows that rule 2 describing the evolution of $X_t$ implies that the norm of the projection of $X_t$ on the hyperplane $H_i$ is bounded by $k^{0.7}$ until time $Bk$.

LEMMA 3.3. *Let $B$ be given by Lemma 3.2 and $j \in \{1, 2, \ldots, i-1, i+1, \ldots, d\}$. Then, for $k$ sufficiently large,*

$$P(\pi_j(X_t) \notin [-k^{0.7}, k^{0.7}] \text{ for some } t \leq Bk) \leq C_5 \exp(-\alpha_5 k^{0.1})$$

*for suitable $C_5 < \infty$ and $\alpha_5 > 0$.*

PROOF. Let $j \neq i$. We observe that, when site $x = X_t + e_j$ is empty, the process $X_t$ jumps to $x$ whenever the particle at $X_t$ gives birth to a particle which is sent to $x$ (type 3 arrow in Figure 5), which occurs at rate $\phi_{ab}$. When site $x$ is occupied, the process $X_t$ jumps to $x$ whenever the particle at $x$ kills the particle at $X_t$ (type 4 arrow in Figure 5), which also occurs at rate $\phi_{ab}$. The same holds by replacing $X_t + e_j$ by $X_t - e_j$. This implies that $\pi_j(X_t)$ is a continuous-time symmetric random walk starting at $\pi_j(X_0) = 0$ run at rate $\phi_{ab}$. This is the main ingredient to prove the lemma. Let $M_k$ denote the number of times $\pi_j(X_t)$ jumps by time $Bk$, and $Y_n$ the discrete-time version of $\pi_j(X_t)$ so that $\pi_j(X_{Bk}) = Y_{M_k}$. The expectation of $M_k$ is given by $m_k = 2Bk\phi_{ab}$. By decomposing the event to be estimated according to whether $M_k > 2m_k$ or $M_k \leq 2m_k$, we obtain

$P(\pi_j(X_t) \notin [-k^{0.7}, k^{0.7}]$ for some $t \leq Bk)$

$\leq P(M_k > 2m_k) + P(Y_n \notin [-k^{0.7}, k^{0.7}]$ for some $n \leq 2m_k; M_k \leq 2m_k)$.

Large deviation estimates imply that the first term can be bounded as follows:

$$P(M_k > 2m_k) \leq C_6 \exp(-\alpha_6 k)$$

for suitable $C_6 < \infty$ and $\alpha_6 > 0$. The second term can be estimated by first using the reflection principle and then Chebyshev's inequality. Given $\theta > 0$, we obtain

$$P(Y_n \notin [-k^{0.7}, k^{0.7}] \text{ for some } n \leq 2m_k; M_k \leq 2m_k)$$
$$\leq 2P(Y_{2m_k} \notin [-k^{0.7}, k^{0.7}]) \leq 2\exp(-\theta k^{0.7})\mathbb{E}\exp(\theta Y_{2m_k})$$
$$\leq 2\exp(-\theta k^{0.7}) \prod_{n=1}^{2m_k} \mathbb{E}\exp(\theta(Y_n - Y_{n-1}))$$
$$\leq 2\exp(-\theta k^{0.7} + 2m_k \log \phi(\theta)),$$

where $\phi(\theta)$ is the moment generating function of $Y_1$. Since $\mathbb{E}Y_1 = 0$ and $\text{Var}\,Y_1 < \infty$, we have that $\log \phi(\theta) \leq C_7 \theta^2$ for some $C_7 < \infty$ and for $\theta$ small enough. By taking $\theta = k^{-0.6}$ in the last expression, we can conclude that

$$P(Y_n \notin [-k^{0.7}, k^{0.7}] \text{ for some } n \leq 2m_k; M_k \leq 2m_k)$$
$$\leq 2\exp(-k^{0.1} + 4C_7\phi_{ab}Bk^{-0.2}) \leq 2\exp(-k^{0.1}/2)$$

for all $k$ sufficiently large. This completes the proof. □

To deduce Lemma 3.1, we first observe that Lemmas 3.2 and 3.3 imply that

$P(\tau_{\omega,i,k} > Bk) \leq P(\sigma_{\omega,i} > Bk)$
$\qquad + P(\pi_j(X_t) \notin [-k^{0.7}, k^{0.7}]$ for some $t \leq Bk$ and some $j \neq i)$
$\qquad \leq C_2 \exp(-\alpha_2 k) + (d-1) \times C_5 \exp(-\alpha_5 k^{0.1})$.

By Lemmas 3.2 and 3.3, we also have

$P(\eta_t(x) = 0$ for all $x \in E_{i,k}$ and some $t \leq \tau_{\omega,i,k})$
$\qquad \leq P(\pi_i(X_t) \leq -k$ for some $t \leq \sigma_{\omega,i})$
$\qquad\quad + P(\pi_j(X_t) \notin [-k^{0.7}, k^{0.7}]$ for some $t \leq \sigma_{\omega,i}$ and some $j \neq i)$
$\qquad \leq P(\pi_i(X_t) \leq -k$ for some $t \leq \sigma_{\omega,i}) + P(\sigma_{\omega,i} > Bk)$
$\qquad\quad + P(\pi_j(X_t) \notin [-k^{0.7}, k^{0.7}]$ for some $t \leq Bk$ and some $j \neq i)$
$\qquad \leq 2C_2 \exp(-\alpha_2 k) + (d-1)C_5 \exp(-\alpha_5 k^{0.1})$

which concludes the proof of Lemma 3.1. □

Note that rule 1 describing the evolution of $X_t$ allows to give the process $X_t$ a drift in the direction of the vector $e_i$ whereas rule 2 implies that the projection of $X_t$ on the orthogonal hyperplane $H_i$ evolves according to a random walk without drift. However, rule 2 can be defined in such a way that the tagged particle has an additional drift towards the $i$th axis. This approach would probably strengthen Lemma 3.1 as follows: a particle at site 0 can be moved of a distance $k$ along the $i$th axis while staying within a distance independent of $k$ of the axis. This result is obviously more difficult to prove than our Lemma 3.1 but would probably shorten the proofs of Lemmas 3.4 and 3.5 below.

The next step is to show that, for given $\varepsilon > 0$, the particle system, when viewed on suitable length and time scales, dominates (in a sense to be made precise) the set of wet sites of an oriented percolation process on the lattice

$$\mathcal{G} = \{(z,n) \in \mathbb{Z}^2 : z + n \text{ is even and } n \geq 0\}$$

in which sites are open with probability $1 - \varepsilon$ [see Durrett (1984), for more details about oriented percolation]. For $\kappa = 1, 2, \ldots$, let $J_\kappa = [-N/\kappa, N/\kappa)^d$ where $N$ is a large integer to be fixed later and let $T > 0$. A site $(z, n) \in \mathcal{G}$ is said to be occupied if there is at least one particle in the translated cube $Nze_1 + J_4$ at time $nT$, that is,

$$\bar{\eta}_{nT}(x) = 1 \qquad \text{for some site } x \in Nze_1 + J_4.$$

The objective is to show that if site $(z, n) \in \mathcal{G}$ is occupied then $(z-1, n+1)$ and $(z+1, n+1)$ are occupied with probability at least $1 - \varepsilon$. Since the evolution rules of the particle system are translation invariant in space and time, it suffices to prove the result when $z = n = 0$. To prepare for the rescaling argument, we need the following lemma, which will be applied repeatedly and in which we prove that, with probability close to 1 when $N$ is large, a particle in $J_{4\kappa}$ can be moved to the smaller cube $J_{16\kappa}$ in less than $dBN$ units of time and before stepping out of $J_\kappa$.

LEMMA 3.4. *Let $\kappa$ be a positive integer and consider the stopping times*

$$\upsilon_{in} = \inf\{t \geq 0 : \bar{\eta}_t(x) = 1 \text{ for some } x \in J_{16\kappa}\} \quad \text{and}$$
$$\upsilon_{out} = \inf\{t \geq 0 : \bar{\eta}_t(x) = 0 \text{ for all } x \in J_\kappa\}.$$

*If there is at least one particle in $J_{4\kappa}$ at time 0, then*

$$P(\upsilon_{in} \leq dBN; \upsilon_{in} < \upsilon_{out}) \geq 1 - C_8 \exp(-\alpha_8 N^{0.1}),$$

*where $B$ is given by Lemma 3.1 and for suitable $C_8 < \infty$ and $\alpha_8 > 0$.*

PROOF. The idea is to apply Lemma 3.1 $d$ times in each direction $i \in \{1, 2, \ldots, d\}$ to move a particle from anywhere in $J_{4\kappa}$ to $J_{16\kappa}$. We refer the reader to Figure 6 for an illustration. To make the argument precise, we let $X \in J_{4\kappa}$ with $\bar{\eta}_0(X) = 1$ and, for $j = 0, 1, \ldots, d$, set

$$D_j = [-jN^{0.7}, jN^{0.7}]^d + \pi_{j+1}(X)e_{j+1} + \pi_{j+2}(X)e_{j+2} + \cdots + \pi_d(X)e_d,$$
$$v_j = \inf\{t \geq 0 : \bar{\eta}_t(x) = 1 \text{ for some } x \in D_j\}$$

and $\omega_j = -\operatorname{sign}(\pi_j(X))$. Note that the sequence starts at $D_0 = X$. By observing that

$$D_d \subset J_{16\kappa} \quad \text{for all } N \text{ sufficiently large} \quad \text{and}$$
$$\varnothing \neq \{x + E_{i,N/2\kappa}\} \cap H_i \subset D_i \quad \text{for all } x \in D_{i-1}$$

and applying Lemma 3.1 with $i = 1, 2, \ldots, d$, $k = N/2\kappa$ and $\omega = \omega_i$, we obtain

$$P(v_{in} > dBN) \leq P(v_d > dBN/2\kappa) \leq \sum_{i=1}^{d} P(v_i - v_{i-1} > BN/2\kappa)$$
$$\leq \sum_{i=1}^{d} P(\tau_{\omega_i, i, k} > BN/2\kappa) \leq C_9 \exp(-\alpha_9 N^{0.1}),$$

where $C_9 = C_1 \times d$ and $\alpha_9 = \alpha_1/(2\kappa)^{0.1}$. Since $x + E_{i,N/2\kappa} \subset J_\kappa$ for all $x \in D_{i-1}$ and the evolution rules of the process are translation invariant in space and time, we can also deduce that

$$P(v_{in} > v_{out}) \leq P(v_d > v_{out}) \leq \sum_{i=1}^{d} P(v_i > v_{out} | v_{i-1} < v_{out})$$
$$\leq \sum_{i=1}^{d} \sup_{x \in D_{i-1}} P(\bar{\eta}_t(z) = 0 \text{ for all } z \in x + E_{i,N/2\kappa})$$

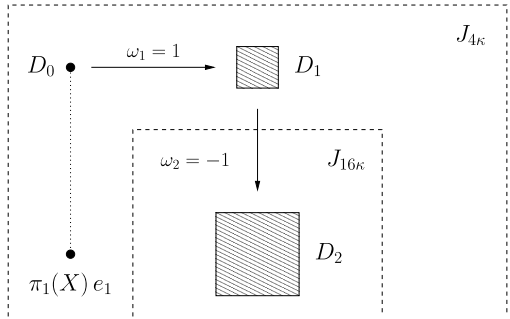

FIG. 6. *Picture of the block argument.*

$$\text{and some } t < v_i \,|\, \bar\eta_0(x) = 1)$$

$$\leq \sum_{i=1}^{d} P(\bar\eta_t(z) = 0 \text{ for all } z \in E_{i,N/2\kappa}$$

$$\text{and some } t < \tau_{\omega_i,i,k} \,|\, \bar\eta_0(0) = 1)$$

$$\leq C_9 \exp(-\alpha_9 N^{0.1}),$$

where $C_9 < \infty$ and $\alpha_9 > 0$ are defined as above. This completes the proof. □

The next lemma provides the main estimate for the rescaling argument.

LEMMA 3.5. *Let $T = (3d+1)BN$ and assume that $\bar\eta_0(X) = 1$ for some $X \in J_4$. Then*

$$P(\bar\eta_T(x) = 0 \text{ for all } x \in J_4 + Ne_1) \leq C_{10}\exp(-\alpha_{10}N^{0.1})$$

*for suitable $C_{10} < \infty$ and $\alpha_{10} > 0$.*

PROOF. For any integer $i \geq 1$, we introduce the stopping times

$$T_j = \begin{cases} \inf\{t > T_{j-1} : \bar\eta_t(x) = 1 \text{ for some } x \in J_{64} + Ne_1\}, & \text{for } j \text{ odd,} \\ \inf\{t > T_{j-1} : \bar\eta_t(x) = 0 \text{ for all } x \in J_{17} + Ne_1\}, & \text{for } j \text{ even,} \end{cases}$$

with the convention $T_0 = 0$. In words, the stopping times $T_j$'s tell us alternately when a particle appears in $J_{64} + Ne_1$ and when $J_{17} + Ne_1$ becomes empty. Let $T^\star$ be the first exit time

$$T^\star = \inf\{t \geq T_1 : \bar\eta_t(x) = 0 \text{ for all } x \in J_4 + Ne_1\}.$$

Then, to establish the lemma, it suffices to prove that $T_1 < T$ and $T^\star > T$ with probability arbitrarily close to 1 when $N$ is large since this indicates that a particle appears in $J_{64} + Ne_1$ by time $T$ and then stays inside $J_4 + Ne_1$ until time $T$.

To estimate the event that $T_1 < T$, we start by applying successively Lemma 3.4 to bring a particle from $J_4$ to $J_{16}$, then from $J_{16}$ to $J_{64}$, and finally from $J_{64}$ to $J_{256}$ in less than $3dBN$ units of time. More precisely, we have

$$P(\bar\eta_t(x) = 0 \text{ for all } x \in J_{256} \text{ and all } t \leq 3dBN) \leq 3 \times C_8\exp(-\alpha_8 N^{0.1}).$$

Once there is a particle in $J_{256}$, we apply Lemma 3.1 with $i = 1$ and $k = N$ to conclude that

$$\begin{aligned}
P(T_1 > T) &= P(T_1 > (3d+1)BN) \\
&\leq 3 \times C_8\exp(-\alpha_8 N^{0.1}) + C_1\exp(-\alpha_1 N^{0.1}) \\
&\leq C_{11}\exp(-\alpha_{11}N^{0.1})
\end{aligned} \tag{9}$$

for $N$ large enough so that $J_{256} + [-N^{0.7}, N^{0.7}]^d \subset J_{64}$ and suitable $C_{11} < \infty$ and $\alpha_{11} > 0$.

We now estimate the event that $T^\star > T$, which is the event that box $J_4 + Ne_1$ contains at least one particle from $T_1$ to $T$. Due to nearest neighbor interactions, there is at least one particle within distance 1 from the boundary of $J_{17} + Ne_1$, and so contained in $J_{16} + Ne_1$, at time $T_{j-1}$ when $j$ is odd. In particular, by applying Lemma 3.4 with $\kappa = 4$, we get

$$(10) \qquad P(T_j > T^\star | T_{j-1} < T^\star) \leq C_8 \exp(-\alpha_8 N^{0.1}) \qquad \text{for } j \text{ odd.}$$

Now, the definition of the $T_j$'s implies that, between time $T_{j-1}$ and time $T_j$ for $j$ even, there is at least one particle in the box $J_{17} + Ne_1$. It follows that

$$(11) \qquad P(T_j > T^\star | T_{j-1} < T^\star) = 0 \qquad \text{for } j \text{ even.}$$

Roughly speaking, (10) and (11) tell us that the number of times a particle appears in $J_{64} + Ne_1$ and $J_{17} + Ne_1$ becomes empty before $J_4 + Ne_1$ becomes empty can be made arbitrarily large by choosing $N$ large. We still need to prove that the number of steps up to time $T$ is smaller than a constant uniform in $N$ by estimating

$$K = \inf\{j \geq 1 : T_j > T\}$$

which is about twice the number of times $J_{17} + Ne_1$ becomes empty by time $T$. Since particles die at rate at most $2d\phi_{ab}$ and the number of deaths required to empty $J_{17} + Ne_1$ starting with at least one particle in $J_{64} + Ne_1$ is larger than $N/32$, we obtain

$$\mathbb{E}[T_j - T_{j-1}] \geq N \times (64 d\phi_{ab})^{-1} \qquad \text{for } j \text{ even.}$$

Using the fact that $T$ is equal to a constant times $N$, large deviation estimates imply the existence of some $m = m(d, \phi_{ab})$ independent of $N$ such that

$$(12) \qquad P(K > m) \leq C_{12} \exp(-\alpha_{12} N)$$

for suitable $C_{12} < \infty$ and $\alpha_{12} > 0$.

Combining estimates (9)–(12), we can conclude that

$$P(\bar{\eta}_T(x) = 0 \text{ for all } x \in J_4 + Ne_1)$$
$$\leq P(T_1 > T) + P(T^\star < T)$$
$$\leq P(T_1 > T) + \sum_{j=2}^{K} P(T_j > T^\star | T_{j-1} < T^\star)$$
$$\leq P(T_1 > T) + P(K > m) + \sum_{j=2}^{m} P(T_j > T^\star | T_{j-1} < T^\star)$$
$$\leq C_{11} \exp(-\alpha_{11} N^{0.1}) + C_{12} \exp(-\alpha_{12} N) + m \times C_8 \exp(-\alpha_8 N^{0.1}).$$

This completes the proof of the lemma. $\square$

Lemma 3.5 implies that, for any $\varepsilon > 0$, the parameters $N$ and $T$ can be chosen in such a way that, provided that site $(0,0)$ is occupied, $(1,1)$ is occupied with probability at least $1 - \varepsilon$. Now, since the range of the interactions is finite, the events described in Lemmas 3.2–3.5 have a finite range of dependence. In particular, with probability at least $1 - \varepsilon$ for $N$ sufficiently large, the configuration of the process at time $T$ in the box $Ne_1 + [-N, N]^d$ only depends upon the graphical representation in the space–time region

$$B_{0,0} = (Ne_1, T) + \{[-2N, 2N]^d \times [0, T]\}.$$

This implies the existence of an event $G_{0,0}$ that has probability at least $1 - 2\varepsilon$, which is measurable with respect to the graphical representation of the process restricted to $B_{0,0}$, and such that, on the event that site $(0,0)$ is occupied and $G_{0,0}$ occurs, site $(1,1)$ is occupied. Finally, since the evolution rules of the process are translation invariant in space and time, we conclude that there exists a collection of events $\{G_{z,n} : (z, n) \in \mathcal{G}\}$ such that, for any $\varepsilon > 0$, the parameters $N$ and $T$ can be chosen in such a way that:

1. The event $G_{z,n}$ is measurable with respect to the graphical representation of the process restricted to the space–time region $B_{z,n} = (Nze_1, nT) + \{[-2N, 2N]^d \times [0, T]\}$.
2. $P(G_{z,n}) \geq 1 - 4\varepsilon$.
3. If $(z, n)$ is occupied and $G_{z,n}$ occurs, then $(z - 1, n + 1)$ and $(z + 1, n + 1)$ are occupied.

This implies that the set of occupied sites dominates the set of wet sites in a 1-dependent oriented site percolation process on $\mathcal{G}$ in which sites are open with probability $p = 1 - 4\varepsilon$ [see Durrett (1995) for a rigorous proof]. This completes the proof of Theorem 2.

To deduce Theorem 1, the last step is to apply a perturbation argument to prove that the previous property holds when $\phi_{aa}$ and $\phi_{bb}$ are positive but small. The probability that a gene originated from a homozygote is sent to the space–time box $B_{z,n}$ can be bounded by

$$(13) \qquad \sum_{x \in [-2N, 2N]^d} \int_0^T \phi \exp(-\phi s) \, ds = (4N + 1)^d \times (1 - \exp(-\phi T)),$$

where $\phi = 2(\phi_{aa} + \phi_{bb})$. The parameters $N$ and $T$ being fixed so that, when $\phi_{aa} = \phi_{bb} = 0$, the set of occupied sites dominates the set of wet sites in a 1-dependent percolation process with parameter $1 - 4\varepsilon$, the rates $\phi_{aa}$ and $\phi_{bb}$ can be chosen so small that the probability in (13) is smaller than $\varepsilon$. With this choice, the set of occupied sites now dominates the set of wet sites in a percolation process with parameter $1 - 5\varepsilon$. To conclude, we take $\varepsilon > 0$ sufficiently small so that percolation occurs, and start the particle system from the "all $ab$" configuration. Then the limit $\nu$ is a stationary distribution. Moreover, the process being coupled with a supercritical oriented percolation process, the density of heterozygotes under the measure $\nu$ is strictly positive. This completes the proof of Theorem 1 when $\phi_{ab} = \phi_{ba}$.

To conclude this section, we now prove Theorem 1 when $\phi_{ab} \neq \phi_{ba}$ and $d = 1$. Since the assumption $\phi_{ab} = \phi_{ba}$ has only been used in the proof of Lemma 3.1, it suffices to extend Lemma 3.1 (see Lemma 3.7 below) to the 1-dimensional case when $\phi_{ab} < 2\phi_{ba}$ and $\phi_{ba} < 2\phi_{ab}$. The difficulty arising from the fact that $\phi_{ab} \neq \phi_{ba}$ is that births of heterozygous individuals onto sites in state $aa$ occur at rate $\phi_{ab}$ whereas births of heterozygous individuals onto sites in state $bb$ occur at rate $\phi_{ba}$. Without loss of generality, we can assume from now on that $\phi_{ab} \geq \phi_{ba}$. To facilitate the writing of the proof, we will use as above the process $\bar\eta_t$ defined in (7). Note however that, due to the general assumption $\phi_{ab} \neq \phi_{ba}$, the process $\bar\eta_t$ is no longer a Markov process. It should be thought of as an auxiliary process that allows us to keep track of heterozygous sites. Let

$$X_t^+ = \sup\{x \in \mathbb{Z} : \bar\eta_t(x) = 1\} \quad \text{and} \quad X_t^- = \inf\{x \in \mathbb{Z} : \bar\eta_t(x) = 1\}$$

denote the right and left edges of the process at time $t$. Then the following lemma holds.

LEMMA 3.6. *Assume that $\phi_{ab} < 2\phi_{ba}$ and start the process with a single particle at site 0. Then there exists a constant $C_{13} > 0$ depending only on the rates $\phi_{ab}$ and $\phi_{ba}$ such that*

$$\liminf_{t \to \infty} t^{-1} X_t^+ > C_{13} \quad and \quad \liminf_{t \to \infty} t^{-1} X_t^- < -C_{13}.$$

PROOF. By symmetry of the evolution rules, it suffices to prove the lemma for $X_t^+$ only. The proof is similar to the proof of Lemma 3.2. We introduce the gap process

$$G_t = \inf\{j \geq 0 : \bar\eta_t(X_t^+ - j - 1) = 1\},$$

TABLE 1
*Evolution rules of $G_t$. In our table, $I_k = \{k, k+1, \ldots\}$, $\circ =$ state 0, and $\bullet =$ state 1. ($0 \to I_1$ indicates a jump from state 0 to some state in $I_1$, etc.) The rows refer to the cases when $G_t = 0$, $G_t = 1$ and $G_t \geq 2$, respectively*

| Configurations | | Transition | Rate | Drift $(\mathbf{X_t^+})$ |
|---|---|---|---|---|
| ∘●●∘∘∘⋯ | ●●●∘∘∘⋯ | $0 \to I_1$ | $\geq \phi_{ab} + \phi_{ba}$ | $\geq (\phi_{ba} - \phi_{ab})/2$ |
| ∘●∘●∘∘⋯ | ●●∘●∘∘⋯ | $1 \to 0$ | $\leq 3\phi_{ab}$ | $\geq \phi_{ba}$ |
| ●∘∘●∘∘⋯ | ●∘∘∘●∘⋯ | $I_2 \to 0$ | $\leq 2\phi_{ab}$ | $\geq \phi_{ba}$ |

and observe that $X_t^+$ jumps to the left at rate $(\phi_{ab}+\phi_{ba})/2$ when $G_t = 0$, and jumps to the right at rate at least $\phi_{ba}$ regardless of the value of $G_t$. Moreover, a straightforward calculation shows that the process $G_t$ makes transitions as indicated in Table 1, where •'s refer to 1's and ∘'s to 0's, and where $I_k$ denotes the set $\{k, k+1, \ldots\}$. The transition from 0 to $I_1$ is obtained by killing one of the particles at $X_t^+$ and $X_t^+ - 1$ (first configuration) or by killing the particle at $X_t^+ - 1$ (second configuration). Both events occur at rate twice the rate at which a particle is killed by one of its neighbors. The last two transitions result from a birth at $X_t^+ - 1$ or $X_t^+ + 1$. Since $\phi_{ab} \geq \phi_{ba}$, upper bounds are reached when each empty site is in state $bb$. These rules imply that

$$P(G_{t+h} = 0) \leq P(G_t = 0) \times [1 - (\phi_{ab}+\phi_{ba})h] + P(G_t \neq 0) \times 3\phi_{ab}h + O(h^2).$$

The same argument as in Lemma 3.2 leads to

$$\liminf_{t \to \infty} t^{-1} X_t^+ \geq \phi_{ba} - \frac{\phi_{ab}+\phi_{ba}}{2} \times \lim_{t \to \infty} P(G_t = 0)$$

$$\geq \phi_{ba} - (\phi_{ab}+\phi_{ba}) \times \frac{3\phi_{ab}}{2(4\phi_{ab}+\phi_{ba})}$$

$$= \frac{2\phi_{ba}^2 + 5\phi_{ab}\phi_{ba} - 3\phi_{ab}^2}{2(4\phi_{ab}+\phi_{ba})}.$$

The roots of $Q(X) = 2X^2 + 5\phi_{ab}X - 3\phi_{ab}^2$ are $-3\phi_{ab}$ and $\phi_{ab}/2$. This implies that the last term of the previous inequality is positive if and only if $\phi_{ab} < 2\phi_{ba}$ and proves the lemma. □

We can now use Lemma 3.6 to demonstrate that with high probability 1's spread out linearly in time and can be found within this growing cone. The second part of Theorem 1 then follows as before from a rescaling argument and we omit the details.

LEMMA 3.7. *Assume that $\phi_{ab} < 2\phi_{ba}$ and $\bar{\eta}_0(0) = 1$. For $\omega \in \{-1, +1\}$, we let*

$$\tau_\omega = \inf\{t \geq 0 : \bar{\eta}_t(\omega k) = 1\}.$$

*Then there exist constants $B < \infty$, $C_{14} < \infty$ and $\alpha_{14} > 0$ such that*

$$P(\tau_\omega > Bk) \leq C_{14} \exp(-\alpha_{14} k) \quad \text{and}$$

$$P(\bar{\eta}_t(x) = 0 \text{ for all } x \in [-k, k] \text{ and some } t \leq \tau_\omega) \leq C_{14} \exp(-\alpha_{14} k).$$

PROOF. By symmetry, it suffices to prove the result for $\omega = +1$. Let

$$Y_t^+ = \sup\{x \leq k : \bar{\eta}_t(x) = 1\} \quad \text{and} \quad Y_t^- = \inf\{x > k : \bar{\eta}_t(x) = 1\}.$$

Since $Y_t^+$ and $Y_t^-$ do not interfere by time $\tau_\omega$, the process $Y_t^+ - Y_0^+$ is equal in distribution to the right edge $X_t^+$ of the particle system starting with a single particle at site 0 by time $\tau_\omega$. In particular, using large deviation estimates, the fact that $Y_0^+ \geq 0$ and Lemma 3.6, we have

$$P(\tau_\omega > Bk) \leq P(X_t^+ < k \text{ for all } t \leq Bk)$$
$$\leq P(X_{Bk}^+ < k) \leq C_{15} \exp(-\alpha_{15} k)$$

for $B = 2/C_{13}$ and suitable $C_{15} < \infty$ and $\alpha_{15} > 0$. Lemma 3.6 also implies that

$$P(\bar\eta_t(x) = 0 \text{ for all } x \in [-k, k] \text{ and some } t \leq \tau_\omega)$$
$$\leq P(X_t^+ < -k \text{ for some } t \geq 0) \leq C_{16} \exp(-\alpha_{16} k).$$

This completes the proof of Lemma 3.7. $\square$

**4. Proof of Theorem 3.** To prove Theorem 3, it is more convenient to consider our particle system as a set of genes evolving on $\mathbb{Z}^d \times \{0,1\}$ rather than a set of genotypes evolving on the lattice $\mathbb{Z}^d$. More precisely, we introduce the Markov process $\xi_t : \mathbb{Z}^d \times \{0,1\} \longrightarrow \{a,b\}$ whose evolution at site $(x, g)$ is given by the following transition rates, where we have set $h = g + 1 \bmod 2$:

$$a \to b \quad \text{at rate } 2\phi_{bb} \sum_{\|x-z\|=1} \mathbb{1}\{\xi(z,g) = \xi(z,h) = b\}$$
$$+ \phi_{ba} \sum_{\|x-z\|=1} \mathbb{1}\{\xi(z,g) = a, \xi(z,h) = b\}$$
$$+ \phi_{ba} \sum_{\|x-z\|=1} \mathbb{1}\{\xi(z,g) = b, \xi(z,h) = a\},$$
$$b \to a \quad \text{at rate } 2\phi_{aa} \sum_{\|x-z\|=1} \mathbb{1}\{\xi(z,g) = \xi(z,j) = a\}$$
$$+ \phi_{ab} \sum_{\|x-z\|=1} \mathbb{1}\{\xi(z,g) = a, \xi(z,h) = b\}$$
$$+ \phi_{ab} \sum_{\|x-z\|=1} \mathbb{1}\{\xi(z,g) = b, \xi(z,h) = a\}.$$

The process $\eta_t$ is a genotype-based model whereas the process $\xi_t$ is a gene-based model, however both describe the same genetic system in the sense that they can be defined on the same probability space in such a way that, for any site $x \in \mathbb{Z}^d$ and any time $t > 0$,

(14) $\quad \eta_t(x) = \xi_t(x,0)\xi_t(x,1) \quad \text{or} \quad \eta_t(x) = \xi_t(x,1)\xi_t(x,0)$

provided the property holds at time 0. In view of (14), it suffices to prove our results for either one of the two processes. We will prove Theorem 3 for the gene-based model.

The proof of Theorem 3 relies on duality techniques supplemented with a standard coupling argument. Let $\phi_a = \min(\phi_{aa}, \phi_{ab})$ and $\phi_b = \max(\phi_{ba}, \phi_{bb})$. The first step is to construct the gene-based model and the biased voter model with parameters $\phi_a$ and $\phi_b$ on the same probability space in such a way that $\xi_t$ has more genes of type $a$ than the biased voter model (denoted by $\zeta_t$ later). To do this, we will construct both processes from the same collections of independent Poisson processes, which is referred to as Harris' graphical representation [Harris (1972)]. We will conclude by showing that $\zeta_t$ converges to the "all $a$" configuration provided $\phi_a > \phi_b$.

To construct the gene-based model $\xi_t$ we put, for all $g, h \in \{0, 1\}$ and all sites $x, z \in \mathbb{Z}^d$ such that $\|x - z\| = 1$, a specific arrow from site $(x, g)$ to site $(z, h)$ at the arrival times of independent Poisson processes with suitable rates. Specifically, for $i, j \in \{a, b\}$, we put, at the arrival times of a Poisson process with parameter $\phi_{ij}$, a type $ij$ arrow from site $(x, g)$ to site $(z, h)$ to indicate that if site $(x, g)$ is occupied by a gene of type $i$ and site $(x, g + 1 \bmod 2)$ is occupied by a gene of type $j$ then site $(z, h)$ becomes occupied by a gene of type $i$. This graphical representation allows us to construct the gene-based model starting from any initial configuration $\xi_0$. Figure 7 shows a realization of the graphical representation of the gene-based process.

To construct the biased voter model $\zeta_t$, we put, for all $g, h \in \{0, 1\}$ and all sites $x, z \in \mathbb{Z}^d$ such that $\|x - z\| = 1$, a type $i$ arrow from site $(x, g)$ to site $(z, h)$ at the arrival times of independent Poisson processes with rate $\phi_i$, $i \in \{a, b\}$. This indicates that if site $(x, g)$ is occupied by a gene of type $i$ then [regardless of the gene present at site $(x, g + 1 \bmod 2)$] site $(z, h)$ becomes occupied by a gene of type $i$.

To construct both processes on the same probability space, we use Poisson processes with rates depending on the signs of $\phi_{aa} - \phi_{ab}$ and $\phi_{bb} - \phi_{ba}$. Table 2 gives an explicit description of the graphical representation in the four possible cases. The type $i$ arrows, $i \in \{a, b\}$, in the construction of the gene-based model have the same interpretation as the ones in the construction of the biased voter model. In particular, a type $i$ arrow in the construction of the gene-based model can be seen as the superposition of a type $ia$ arrow and a type $ib$ arrow. In view of the well-known properties of the Poisson process, each of these graphical representations produces the desired flip rates to construct the processes we are interested in. Moreover, it follows directly from Table 2 that

(15) $\quad \zeta_t(x, i) = a \implies \xi_t(x, i) = a \quad$ for all $x \in \mathbb{Z}^d$ and $i \in \{0, 1\}$

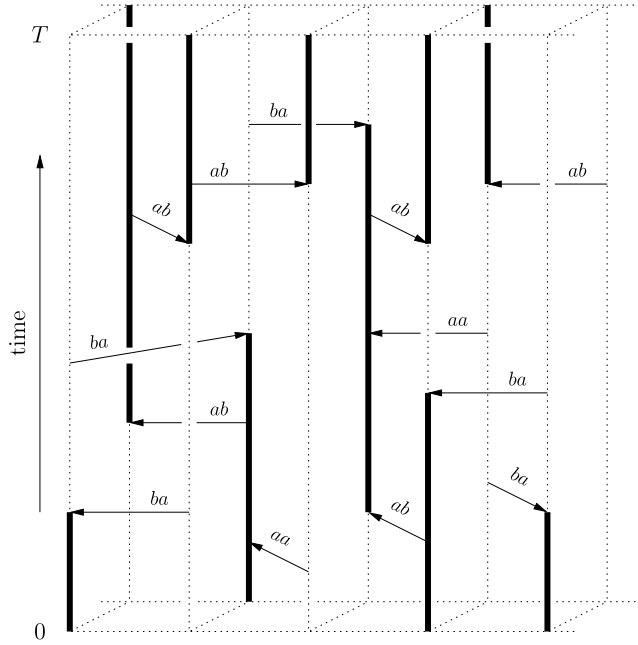

Fig. 7. *Construction of the gene-based model from the Harris' graphical representation. The thick lines refer to sites occupied by a gene of type a.*

at any time $t > 0$ provided this holds at time 0. Note, however, that in our coupling the Poisson processes with rate $\phi_{ij}$, $i,j \in \{a,b\}$, are no longer independent.

Assume that $\phi_a > \phi_b$ and that infinitely many genes of type $a$ are present at time 0. The last step is to prove that the biased voter model $\zeta_t : \mathbb{Z}^d \times \{0,1\} \longrightarrow \{a,b\}$ with parameters $\phi_a$ and $\phi_b$ converges almost surely to the "all $a$" configuration. The proof relies on duality techniques. We construct the process by drawing a type $a$ arrow (respectively, an unlabeled arrow) from site $(x,g)$ to site $(z,h)$ at the arrival times of a Poisson process with rate $\phi_a - \phi_b$ (respectively, $\phi_b$). Unlabeled arrows indicate that the gene at site $(z,h)$ is replaced by the gene at site $(x,g)$ regardless of the type of both genes. The dual process is constructed by going backwards in time and jumps at the tips of unlabeled arrows while branches at the tips of type $a$ arrows as indicated in Figure 8. More precisely, if $\hat{\zeta}_s$ denotes the dual process at dual time $s$ then:

1. if there is an unlabeled arrow from $(x,g)$ to $(z,h)$ at time $t-s$ for some $(z,h) \in \hat{\zeta}_{s-}$, then $\hat{\zeta}_s$ can be deduced from $\hat{\zeta}_{s-}$ by removing site $(z,h)$ and adding site $(x,g)$ while

Table 2
Couplings of the gene-based model and the biased voter model

| Rate | Gene-based | Biased voter |
|---|---|---|
| | $\phi_{aa} \leq \phi_{ab}$ and $\phi_{bb} \leq \phi_{ba}$ | |
| $\phi_{aa}$ | $\xrightarrow{a}$ | $\xrightarrow{a}$ |
| $\phi_{ab} - \phi_{aa}$ | $\xrightarrow{ab}$ | |
| $\phi_{bb}$ | $\xrightarrow{b}$ | $\xrightarrow{b}$ |
| $\phi_{ba} - \phi_{bb}$ | $\xrightarrow{ba}$ | $\xrightarrow{b}$ |
| | $\phi_{aa} \leq \phi_{ab}$ and $\phi_{bb} \geq \phi_{ba}$ | |
| $\phi_{aa}$ | $\xrightarrow{a}$ | $\xrightarrow{a}$ |
| $\phi_{ab} - \phi_{aa}$ | $\xrightarrow{ab}$ | |
| $\phi_{ba}$ | $\xrightarrow{b}$ | $\xrightarrow{b}$ |
| $\phi_{bb} - \phi_{ba}$ | $\xrightarrow{bb}$ | $\xrightarrow{b}$ |
| | $\phi_{aa} \geq \phi_{ab}$ and $\phi_{bb} \leq \phi_{ba}$ | |
| $\phi_{ab}$ | $\xrightarrow{a}$ | $\xrightarrow{a}$ |
| $\phi_{aa} - \phi_{ab}$ | $\xrightarrow{aa}$ | |
| $\phi_{bb}$ | $\xrightarrow{b}$ | $\xrightarrow{b}$ |
| $\phi_{ba} - \phi_{bb}$ | $\xrightarrow{ba}$ | $\xrightarrow{b}$ |
| | $\phi_{aa} \geq \phi_{ab}$ and $\phi_{bb} \geq \phi_{ba}$ | |
| $\phi_{ab}$ | $\xrightarrow{a}$ | $\xrightarrow{a}$ |
| $\phi_{aa} - \phi_{ab}$ | $\xrightarrow{aa}$ | |
| $\phi_{ba}$ | $\xrightarrow{b}$ | $\xrightarrow{b}$ |
| $\phi_{bb} - \phi_{ba}$ | $\xrightarrow{bb}$ | $\xrightarrow{b}$ |

2. if there is a type $a$ arrow from $(x, g)$ to $(z, h)$ at time $t - s$ for some $(z, h) \in \hat{\zeta}_{s^-}$, then $\hat{\zeta}_s$ can be deduced from $\hat{\zeta}_{s^-}$ by adding site $(x, g)$.

The convergence to the "all $a$" configuration follows from the analysis of the dual process. A stronger version of this result is proved in Bramson and Griffeath (1980, 1981) for the biased voter model on $\mathbb{Z}^d$. Their proof, however, easily extends to the process on $\mathbb{Z}^d \times \{0, 1\}$. This, together with (14) and (15), concludes the proof of Theorem 3.

**5. Proof of Theorem 4.** Similar to the proof of Theorem 1, the proof of Theorem 4 relies on a rescaling argument. First, we set $F_N = [-N, N]^d$ for $N \geq 0$ and denote by $P_{F_N}$ the law of the process starting with all sites in $F_N$ occupied by individuals of genotype $aa$, and all sites outside $F_N$ occupied by individuals of genotype $bb$. The first step is to prove the following lemma.

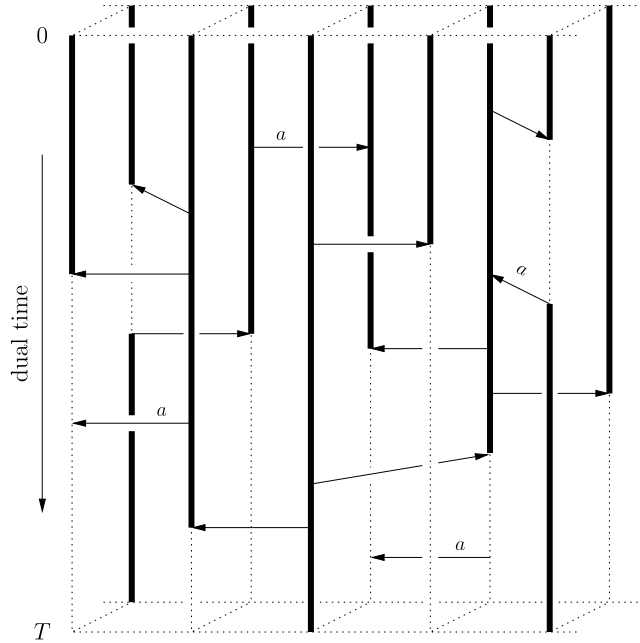

Fig. 8. *Picture of the dual coalescing branching random walks.*

Lemma 5.1. *Assume (4) and $\phi_{ba} = 0$. Then*

$$P_{F_N}(\eta_t(0) \neq aa \text{ at some time } t \geq 0) \leq C_{17} \exp(-\alpha_{17} N)$$

*for suitable constants $C_{17} < \infty$ and $\alpha_{17} > 0$.*

The proof of Lemma 5.1 is made difficult by the lack of a dual process. To estimate the probability that a gene of type $b$ reaches site 0, we will construct a path $\gamma \subset \mathbb{Z}^d$ starting outside the cube $F_N$ and ending in the neighborhood of site 0 along which genes of type $b$ spread. The construction will be done by going backwards in time. By investigating the evolution of the process along the invasion path $\gamma$ by going forward in time, we will prove that the probability that such a path exists is smaller than a constant times $c^{-K}$ for some $c > 2d$, where $K$ denotes the length of $\gamma$. We will conclude the proof of Lemma 5.1 by summing over all the possible paths starting at site $x$ for some $x \in \mathbb{Z}^d \setminus F_N$ and ending in the neighborhood of site 0.

Lemma 5.1 implies that, with probability arbitrarily close to 1 when $N$ is large, the individual at site 0 is of type $aa$ at any time $t \geq 0$. The second step is to use this individual as a source to fill the cube $F_{3N}$ with genes of type $a$. This occurs by some time $T < \infty$ with probability arbitrarily close to 1 for $N$ large. By applying the same arguments as in the proof of Lemma 5.1 to all the individuals in $F_{2N}$ instead of the individual at site 0, we will deduce the following lemma.

LEMMA 5.2. *Assume (4) and $\phi_{ba} = 0$. For any $\varepsilon > 0$,*

$$P_{F_N}(\eta_T(x) \neq aa \text{ for some } x \in F_{2N}) \leq 2\varepsilon$$

*for a suitable integer $N \geq 0$ and a (deterministic) time $T < \infty$.*

Lemma 5.2 implies the existence of a partition of $\mathbb{Z}^d$ into cubes of side length $2N + 1$ such that the set of cubes filled with $a$'s at time $2nT$, $n \geq 0$, dominates the set of wet sites at level $2n$ of a 1-dependent oriented percolation process with parameter $1 - 2\varepsilon$. This can be proved by applying the same arguments as in the proof of Theorem 1 (for more details, see our construction after the proof of Lemma 3.5). To prove that genes of type $a$ actually outcompete genes of type $b$, namely that there is an in-all-directions expanding region which is void of genes of type $b$, we apply a result from Durrett (1992) that shows that sites that are not occupied do not percolate when $\varepsilon$ is close enough to 0. Since genes of either type cannot appear spontaneously, once a region is void of one type, this type can only reappear in the region through invasion from the outside. To complete the proof of Theorem 4, the last step is to apply a perturbation argument similar to the one described in the proof of Theorem 1 to extend the result to the case when $\phi_{ba}$ is small. We now prove Lemmas 5.1 and 5.2 in detail.

5.1. *Construction of the invasion path $\gamma$.* The first step in proving Lemma 5.1 is to demonstrate the existence of a path $\gamma \subset \mathbb{Z}^d$ along which genes of type $b$ spread in the case when the stopping time

$$T_N = \inf\{t \geq 0 : \eta_t(0) \neq aa\}$$

is finite. We first let the process starting with genes of type $aa$ in $F_N$ evolve until time $T_N$ and then do the construction inductively by going backwards in time.

We set $S_0 = T_N$. The definition of $S_0$ implies that there is a gene of type $b$ that is sent to site 0 at time $S_0$. This gene must have originated from a site $\Gamma_0$ with $\|\Gamma_0\| = 1$. Since $\phi_{ba} = 0$, we deduce that site $\Gamma_0$ is necessarily in state $bb$ at time $S_0$. Then we let

$$S_1 = \sup\{s \leq S_0 : \eta_s(\Gamma_0) = aa\}.$$

Again, there is a site $\Gamma_1$ with $\|\Gamma_1 - \Gamma_0\| = 1$ which is in state $bb$ at time $S_1$. Let

$$S_2 = \sup\{s \leq S_1 : \eta_s(\Gamma_1) = aa\}.$$

This procedure allows us to construct a sequence $(\Gamma_i, S_i)$ such that time $S_i$ is almost surely finite and site $\Gamma_i$ is in state $bb$ at time $S_i$. The construction stops at $i = K$ where

$$K = \min\{i \geq 0 : \eta_s(\Gamma_i) \neq aa \text{ for all } s \in [0, S_i]\}.$$

Note that this implies that $\Gamma_K \notin F_N$ and that $K \geq N$. Our construction implies that

$$0 < S_K < S_{K-1} < \cdots < S_0 = T_N.$$

Since it will be more convenient to work forward in time, we construct the invasion path $\gamma$ by reversing the direction of time, that is we set

$$\gamma(i) = \Gamma_{K-i} \quad \text{and} \quad s_i = S_{K-i}$$

for $i = 0, 1, \ldots, K$, and $\gamma = (\gamma(0), \gamma(1), \ldots, \gamma(K))$.

5.2. *Evolution of the genotypes along $\gamma$.* The next step is to define a process $\omega_t$ starting at $\omega_0 = 0$ such that $\gamma(\omega_t)$ follows the invasion of genes of type $b$ in the cube $F_N$. The process $\omega_t$ is stochastically smaller than a process whose probability to reach $K$ decreases exponentially in $K$ (see Lemma 5.3 below). This allows us to conclude that the probability that the path $\gamma$ exists decreases exponentially fast with its length, which is the key to prove Lemma 5.1.

The process $\omega_t$ starts at $\omega_0 = 0$ and has values in $\{0, 1, \ldots, K\}$. The transition rates depend on whether site $\gamma(\omega_t + 1)$ is in state $aa$ or not.

1. If $\gamma(\omega_t + 1)$ is in state $aa$ (at time $t$), then

$$\omega_t \to \begin{cases} \omega_t + 1, & \text{when the state of } \gamma(\omega_t + 1) \text{ jumps to } ab, \\ \sup\{i \leq \omega_t - 1 : \eta_t(\gamma(i)) \neq aa\}, \\ & \text{when the state of } \gamma(\omega_t) \text{ jumps to } aa. \end{cases}$$

2. If $\gamma(\omega_t + 1)$ is not in state $aa$ (at time $t$), then

$$\omega_t \to \sup\{i \leq \omega_t - 1 : \eta_t(\gamma(i)) \neq aa\} \quad \text{when the state of } \gamma(\omega_t) \text{ jumps to } aa.$$

Note that $\omega_t$ is well defined until time $T_N$ and $\sup\{i \leq \omega_t - 1 : \eta_t(\gamma(i)) \neq aa\} \geq 0$ until this stopping time condition on the event that the invasion path $\gamma$ exists. Our construction also implies that $\omega_t$ jumps from $i$ to $i + 1$ at time $s_i$ (see the left-hand side of Figure 9) so that $\omega_t$ reaches state $K$ at time $s_K = T_N$, and that site $\gamma(\omega_t)$ is in state $ab$ or $bb$ up to time $s_K$. To keep track of the genotype, we define a process $\bar{\omega}_t$ with state space $(\mathbb{Z}/2) \cap (-1, K]$ by letting

$$\bar{\omega}_t = \omega_t - \tfrac{1}{2} \mathbb{1}\{\eta_t(\gamma(\omega_t)) = ab, \eta_t(\gamma(\omega_t + 1)) = aa\}.$$

See the right-hand side of Figure 9. We now describe the process $\bar{\omega}_t$ in details. The transition rates of the process depend on the state of sites $\gamma(\omega_t)$ and $\gamma(\omega_t + 1)$, and we distinguish four different cases that we call configurations 1–4 (see Figure 10 for an illustration):

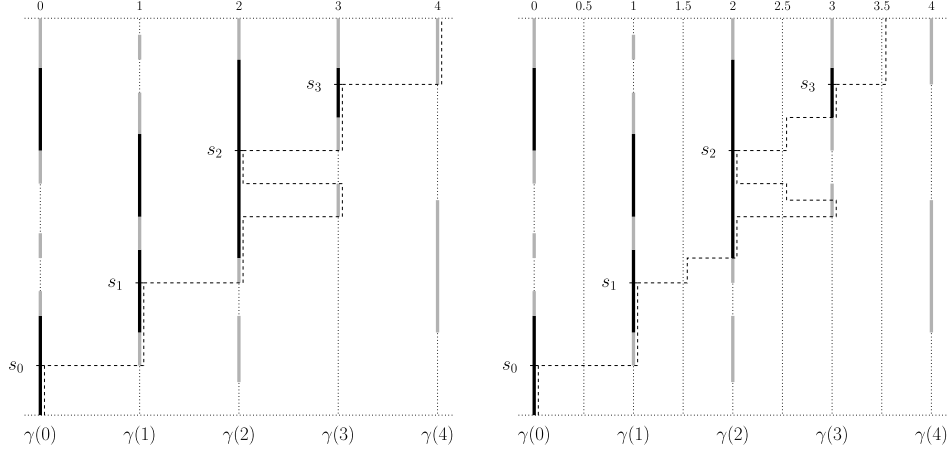

FIG. 9. *Illustration of $\omega_t$ (on the left-hand side) and $\bar{\omega}_t$ (on the right-hand side) in dashed lines from a realization of the particle system. In both pictures, the horizontal axes represents both the state space of $\gamma$ (on the bottom) and the state space of the processes $\omega_t$ and $\bar{\omega}_t$ (on the top). Vertical thick black lines refer to sites in state bb, thick grey lines to sites in state ab and dotted lines to sites in state aa.*

1. Config. 1—site $\gamma(\omega_t)$ is in state *bb* and site $\gamma(\omega_t + 1)$ is in state *aa*.
   In this case, $\bar{\omega}_t$ can either jump to $\bar{\omega}_t + 1/2$ or $\bar{\omega}_t + 1$ when the state of site $\gamma(\omega_t + 1)$ jumps to *ab*, or jump to $\bar{\omega}_t - 1/2$ when the state of site

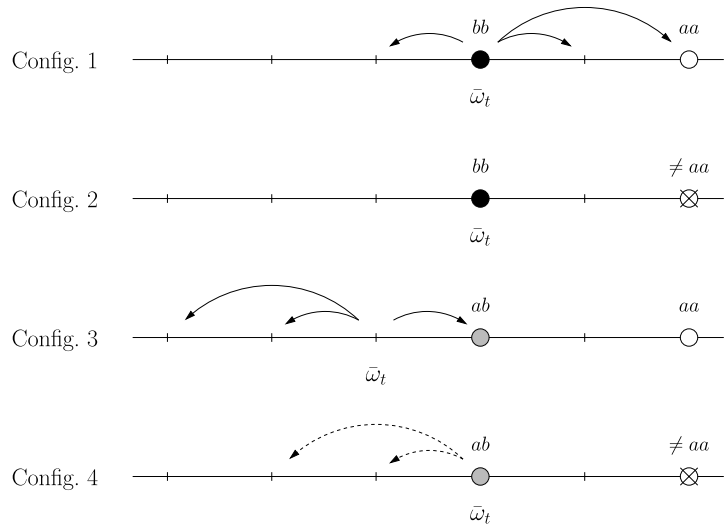

FIG. 10. *Pictures of the configurations 1–4.*

$\gamma(\omega_t)$ jumps to $ab$ which implies that

$$\bar{\omega}_t \to \begin{cases} \leq \bar{\omega}_t + 1, & \text{at rate at most } 4d\phi_{bb}, \\ \bar{\omega}_t - 1/2, & \text{at rate at least } 2\phi_{aa}. \end{cases}$$

2. Config. 2—site $\gamma(\omega_t)$ is in state $bb$ and site $\gamma(\omega_t + 1)$ is not in state $aa$.
   In this case, the process $\bar{\omega}_t$ is frozen.
3. Config. 3—site $\gamma(\omega_t)$ is in state $ab$ and site $\gamma(\omega_t + 1)$ is in state $aa$.
   In this case, $\bar{\omega}_t$ can either jump to $\bar{\omega}_t + 1/2$ when the state of site $\gamma(\omega_t)$ jumps to $bb$ or the state of site $\gamma(\omega_t + 1)$ jumps to $ab$, or jump to the left when the state of site $\gamma(\omega_t)$ jumps to $aa$ which implies that

$$\bar{\omega}_t \to \begin{cases} \bar{\omega}_t + 1/2, & \text{at rate at most } 6d\phi_{bb}, \\ \leq \bar{\omega}_t - 1/2, & \text{at rate at least } \phi_{aa}. \end{cases}$$

4. Config. 4—site $\gamma(\omega_t)$ is in state $ab$ and site $\gamma(\omega_t + 1)$ is not in state $aa$.
   In this case, the process $\bar{\omega}_t$ can only jump to the left.

The process $\bar{\omega}_t$ is well defined and $\geq -1/2$ until time $T_N$ condition on the event that the invasion path $\gamma$ exists. However, when $\omega_t = 0$ and the configuration is of type 3, the state of $\gamma(0)$ may still jump to $aa$ at the rate indicated above, thus contradicting the existence of $\gamma$. This motivates, in order to bound the probability that $\gamma$ exists, the introduction of the Markov process $X_t$ with state space $\mathbb{Z}/2$ and with transition rates

$$X_t \to \begin{cases} \lfloor X_t \rfloor + 1, & \text{at rate } 6d\phi_{bb}, \\ X_t - 1/2, & \text{at rate } \phi_{aa}, \end{cases}$$

where $\lfloor \cdot \rfloor$ denotes the integer part. The choice of the transition rates of $X_t$ and the observation above imply that the probability that the invasion path $\gamma$ exists is bounded by the probability that, starting at 0, the process $X_t$ reaches state $K$ before state $-1$. To estimate the latter, we consider the discrete-time version of $X_t$, which is the Markov process $Y_n$ with transition probabilities

$$(16) \qquad Y_{n+1} = \begin{cases} \lfloor Y_n \rfloor + 1, & \text{with probability } r, \\ Y_n - 1/2, & \text{with probability } l, \end{cases}$$

where the probabilities $r$ (right) and $l$ (left) are given by

$$r = \frac{6d\phi_{bb}}{\phi_{aa} + 6d\phi_{bb}} \quad \text{and} \quad l = \frac{\phi_{aa}}{\phi_{aa} + 6d\phi_{bb}}.$$

Then, denoting by $\tau_i$ the first time the process $Y_n$ hits state $i$, we have

$$(17) \qquad \begin{aligned} & P(\gamma = (\gamma(0), \gamma(1), \ldots, \gamma(K))) \text{ is an invasion path}) \\ & \leq P(\tau_K < \tau_{-1} | Y_0 = 0), \end{aligned}$$

which decreases exponentially in $K$ according to the following lemma.

LEMMA 5.3. *Assume (4). Then there exists a constant $c > 2d$ such that*
$$P(\tau_K < \tau_{-1}|Y_0 = 0) = (c-1) \cdot (c^{K+1} - 1)^{-1} \leq c^{-K}.$$

PROOF. Let $P_i$ denote the law of the process starting at $Y_0 = i$. Then the transition probabilities of the process indicated in (16) imply that
$$P_i(\tau_{i-1} < \tau_{i+1}) = l^2 \sum_{j=0}^{\infty} (l \cdot r)^j = \frac{l^2}{1 - lr} \quad \text{and} \quad P_i(\tau_{i-1} > \tau_{i+1}) = \frac{r}{1 - lr}$$

for any $i \in \mathbb{Z}$. Let $p_i = P_i(\tau_K < \tau_{-1})$ and $w(i) = (l^2/r)^i$. By decomposing according to whether the process first jumps to $i - 1$ or $i + 1$, and using the fact that
$$(l^2 + r)w(i) = l^2 w(i-1) + r w(i+1),$$

we deduce that $w(Y_n)$ is a martingale. The martingale stopping theorem then implies that
$$w(i) = p_i w(K) + (1 - p_i) w(-1)$$

from which it follows that
$$p_i = (w(i) - w(-1)) \cdot (w(K) - w(-1))^{-1} = (c^{i+1} - 1) \cdot (c^{K+1} - 1)^{-1}$$

with $c = l^2/r$. This implies that $P_0(\tau_K < \tau_{-1}) = p_0 = (c - 1) \cdot (c^{K+1} - 1)^{-1}$. Finally, since
$$l^2/r > 2d \quad \text{if and only if} \quad \phi_{aa}^2 > 12d^2 \phi_{bb}(\phi_{aa} + 6d\phi_{bb}),$$
$$\text{if and only if} \quad \phi_{aa} > 6d^2 \left(1 + \sqrt{1 + \frac{2}{d}}\right)\phi_{bb}$$

which is condition (4) in Theorem 4; the lemma follows. □

5.3. *Proof of Lemma 5.1.* First of all, combining Lemma 5.3 with (17), we obtain
$$P(\gamma = (\gamma(0), \gamma(1), \ldots, \gamma(K)) \text{ is an invasion path}) \leq c^{-K}$$

with $c > 2d$ given in Lemma 5.3. Now, let $x \in \partial F_{L+1} = F_{L+1} \setminus F_L$ for some $L \geq 0$ and denote the set of paths starting at $x$ and ending in the neighborhood of 0 by $\Theta_x$. By the triangle inequality, for any $\gamma = (\gamma(0), \gamma(1), \ldots, \gamma(K)) \in \Theta_x$,
$$L \leq \|\gamma(K) - \gamma(0)\| \leq \sum_{i=0}^{K-1} \|\gamma(i+1) - \gamma(i)\| = K,$$

that is, each path in $\Theta_x$ has length at least $L$. Furthermore, since each site has $2d$ nearest neighbors, the number of possible paths of length $K$ is bounded by $(2d)^K$. It follows that

$P(\text{there exists } \gamma \in \Theta_x \text{ such that } \gamma \text{ is an invasion path})$

$$\leq \sum_{K=L}^{\infty} P(\text{there exists } \gamma \in \Theta_x \text{ such that } \gamma \text{ is an invasion path of length } K)$$

$$\leq \sum_{K=L}^{\infty} (2d/c)^K = \frac{c}{c-2d} \times (2d/c)^L.$$

Finally, since the birth of a gene of type $b$ at site $0$ implies the existence of an invasion path starting at some site $x \notin F_N$ for the process starting with all sites in $F_N$ occupied by individuals of genotype $aa$, by summing over all sites $x \notin F_N$, we obtain

$P_{F_N}(T_N < \infty)$

$$\leq \sum_{x \notin F_N} P(\text{there exists } \gamma \in \Theta_x \text{ such that } \gamma \text{ is an invasion path})$$

$$\leq \sum_{L=N}^{\infty} \sum_{x \in \partial F_{L+1}} P(\text{there exists } \gamma \in \Theta_x \text{ such that } \gamma \text{ is an invasion path})$$

$$\leq \sum_{L=N}^{\infty} ((2L+3)^d - (2L+1)^d) \times \frac{c}{c-2d} \times (2d/c)^L$$

$$\leq C_{17}\, (2d/c)^N \leq C_{17} \exp(-\alpha_{17} N)$$

for suitable $C_{17} < \infty$ and $\alpha_{17} > 0$. This completes the proof of Lemma 5.1.

5.4. *Proof of Lemma 5.2.* Lemma 5.2 is a straightforward consequence of Lemma 5.1. Let $\varepsilon > 0$ and denote by $T^\star_{3N}$ the first time the cube $F_{3N}$ is void of genes of type $b$, that is

$$T^\star_{3N} = \inf\{t \geq 0 : \eta_t(x) = aa \text{ for all } x \in F_{3N}\}.$$

By applying Lemma 5.1 to each of the sites belonging to $\partial F_{2N+1}$ at time $T^\star_{3N}$, we obtain

$P(\eta_t(x) \neq aa \text{ for some } x \in F_{2N} \text{ at some time } t \geq T^\star_{3N} \mid T^\star_{3N} < \infty)$

$$\leq P(\eta_t(x) \neq aa \text{ for some } x \in \partial F_{2N+1} \text{ at some time } t \geq T^\star_{3N} \mid T^\star_{3N} < \infty)$$

$$\leq ((2N+3)^d - (2N+1)^d) \times C_{17} \exp(-\alpha_{17} N)$$

$$\leq C_{18} \exp(-\alpha_{18} N)$$

for suitable $C_{18} < \infty$ and $\alpha_{18} > 0$. Pick $N$ sufficiently large so that
$$C_{17} \exp(-\alpha_{17} N) + C_{18} \exp(-\alpha_{18} N) \leq \varepsilon.$$
Now that $N$ is fixed, we will define a time $T$ satisfying Lemma 5.2 as a function of $N$. First of all, we observe that if site 0 is in state $aa$ at some time $t$ then there is a positive probability (uniform in $t$) that the cube $F_{3N}$ is void of genes of type $b$ at time $t+1$. This implies that
$$P(T^\star_{3N} < \infty | T_N = \infty) = 1.$$
Since the sequence of events $\{T^\star_{3N} \leq n\}$, $n \geq 0$, is nondecreasing, we can apply Beppo–Levi theorem to obtain the existence of a large enough deterministic time $T < \infty$ such that
$$P(T^\star_{3N} < T | T_N = \infty) = 1 - \varepsilon.$$
In conclusion,
$$\begin{aligned}
P_{F_N}(\eta_T&(x) \neq aa \text{ for some } x \in F_{2N}) \\
&\leq P_{F_N}(T_N < \infty) + P(\eta_T(x) \neq aa \text{ for some } x \in F_{2N}; T_N = \infty) \\
&\leq P_{F_N}(T_N < \infty) + P(T_N = \infty; T^\star_{3N} \geq T) \\
&\quad + P(\eta_T(x) \neq aa \text{ for some } x \in F_{2N}; T^\star_{3N} < T) \\
&\leq P_{F_N}(T_N < \infty) + P(T_N = \infty; T^\star_{3N} \geq T) \\
&\quad + P(\eta_t(x) \neq aa \text{ for some } x \in F_{2N} \text{ at some time } t \geq T^\star_{3N}; T^\star_{3N} < \infty) \\
&\leq \varepsilon + C_{17} \exp(-\alpha_{17} N) + C_{18} \exp(-\alpha_{18} N) \leq 2\varepsilon.
\end{aligned}$$
This completes the proof of Lemma 5.2 and Theorem 4.

**6. Proof of Theorem 6.** Theorem 6 predicts convergence to the "all $aa$" configuration for a suitable choice of the parameters and for the 1-dimensional process starting from any initial configuration with infinitely many individuals of genotype $aa$. In Sections 6.1 and 6.2, we will prove for the process starting with only $aa$ on the half line $(-\infty, 0]$ that $\lim_{t \to \infty} r_t = \infty$ where
$$r_t = \sup\{x \in \mathbb{Z} : \eta_t(z) = aa \text{ for all } z \leq x\}.$$
This implies that (i) starting with a single individual of genotype $aa$ at site 0
$$\lim_{t \to \infty} Z^-_t = -\infty \quad \text{and} \quad \lim_{t \to \infty} Z^+_t = +\infty,$$
almost surely on the event $\{Z^-_t \leq 0 \leq Z^+_t \text{ for all } t \geq 0\}$ where
$$\begin{aligned}
Z^-_t &= \inf\{x \leq 0 : \eta_t(z) = aa \text{ for all } x \leq z \leq 0\}, \\
Z^+_t &= \sup\{x \geq 0 : \eta_t(z) = aa \text{ for all } x \geq z \geq 0\},
\end{aligned}$$

and that (ii) $P(Z_t^- \leq 0 \leq Z_t^+$ for all $t \geq 0) > C_{19} > 0$. From (i) and (ii), it is straightforward to deduce that, with probability 1, the process starting with infinitely many individuals of genotype $aa$ converges to the "all $aa$" configuration. Coming back to our initial objective, we assume from now on that the process starts with only individuals of genotype $aa$ in the interval $(-\infty, 0]$. To show that $\lim_{t\to\infty} r_t = \infty$, we first set $K_t = s_t - r_t - 1$ where

$$s_t = \inf\{x \geq r_t + 1 : \eta_t(x) = bb\}$$

and observe that the process $K_t$ evolves as indicated in Table 3. The proof is carried out in two steps according to whether $\phi_{aa} > \phi_{ba}$ or $\phi_{aa} \leq \phi_{ba}$.

6.1. $\phi_{aa} > \phi_{ba}$. In this case, the process $r_t$ has a positive drift when $K_t \neq 0$ (see Table 3) and the proof consists in estimating the fraction of time $K_t \neq 0$ when

$$\phi_{aa} > \phi_{bb} + \phi_{ba} + \sqrt{\phi_{bb}\phi_{ba}}.$$

First of all, we observe that, since $\phi_{aa} > \max(\phi_{ba}, \phi_{bb})$, the process $K_t$ jumps from $I_1$ to 0 at rate at most $2\phi_{aa}$, which implies that

$$P(K_{t+h} = 0) \leq P(K_t = 0) \times [1 - 2(\phi_{aa} + \phi_{bb})h] \\ + P(K_t \neq 0) \times 2\phi_{aa}h + O(h^2).$$

By first taking the limit as $h \to 0$ and then the limit as $t \to \infty$, it follows that

$$\frac{d}{dt} P(K_t = 0) \leq 2\phi_{aa} - (4\phi_{aa} + 2\phi_{bb})P(K_t = 0)$$

and then that

$$\limsup_{t\to\infty} P(K_t = 0) \leq \frac{\phi_{aa}}{2\phi_{aa} + \phi_{bb}}.$$

TABLE 3
*Evolution rules of $K_t$. In our table, $I_k = \{k, k+1, \ldots\}$, $\circ =$ state aa, $\bullet =$ state ab and $\bullet =$ state bb*

| Configuration | Transition | Rate | Drift ($r_t$) |
|---|---|---|---|
| ∘∘∘●⋯ | $0 \to I_1$ | $\geq 2(\phi_{aa} + \phi_{bb})$ | $\geq -2\phi_{bb}$ |
| ∘∘●●⋯ | $1 \to 0$ | $\phi_{aa} + \phi_{bb}$ | $\geq \phi_{aa} - \phi_{ba}$ |
| ∘∘●●⋯ | $1 \to I_2$ | $\geq \phi_{ab} + \phi_{ba}$ | $\geq \phi_{aa} - \phi_{ba}$ |
| ∘∘●●⋯ | $I_2 \to 0$ | $\phi_{ba}/2$ | $\geq \phi_{aa} + \phi_{ab}/2 - \phi_{ba}$ |
| ∘∘●∘⋯ | $I_2 \to 0$ | $\leq 2\phi_{aa}$ | $\geq 4\phi_{aa} - \phi_{ba}$ |
| ∘∘●●⋯ | $I_2 \to 1$ | $\leq \phi_{aa} + \phi_{bb} + \phi_{ab}/2 + \phi_{ba}/2$ | $\geq \phi_{aa} + \phi_{ab}/2 - \phi_{ba}$ |
| ∘∘●∘⋯ | $I_2 \to 1$ | $0$ | $\geq 4\phi_{aa} - \phi_{ba}$ |

Since the drift is bounded from below by $\phi_{aa} - \phi_{ba}$ when $K_t \neq 0$ (see Table 3), we obtain

$$\liminf_{t \to \infty} \frac{r_t}{t} \geq -2\phi_{bb} \limsup_{t \to \infty} P(K_t = 0) + (\phi_{aa} - \phi_{ba}) \liminf_{t \to \infty} P(K_t \neq 0)$$

$$\geq -2\phi_{bb} \frac{\phi_{aa}}{2\phi_{aa} + \phi_{bb}} + (\phi_{aa} - \phi_{ba}) \frac{\phi_{aa} + \phi_{bb}}{2\phi_{aa} + \phi_{bb}}$$

$$= [-2\phi_{bb} + (\phi_{aa} - \phi_{ba})(1 + \phi_{bb}/\phi_{aa})] \frac{\phi_{aa}}{2\phi_{aa} + \phi_{bb}}$$

$$= \frac{\phi_{aa}^2 - (\phi_{bb} + \phi_{ba})\phi_{aa} - \phi_{bb}\phi_{ba}}{2\phi_{aa} + \phi_{bb}}.$$

This leads to the requirement that $\phi_{aa}^2 - (\phi_{bb} + \phi_{ba})\phi_{aa} - \phi_{bb}\phi_{ba} \geq 0$, which is equivalent to

$$\phi_{aa} \geq \tfrac{1}{2}(\phi_{bb} + \phi_{ba}) + \sqrt{\phi_{bb}\phi_{ba} + \tfrac{1}{4}(\phi_{bb} + \phi_{ba})^2}.$$

Since, by subadditivity of the function square root,

$$\sqrt{\phi_{bb}\phi_{ba} + \tfrac{1}{4}(\phi_{bb} + \phi_{ba})^2} \leq \sqrt{\phi_{bb}\phi_{ba}} + \tfrac{1}{2}(\phi_{bb} + \phi_{ba}),$$

this shows that condition (5) implies that $\lim_{t \to \infty} r_t = \infty$.

6.2. $\phi_{aa} \leq \phi_{ba}$. The process $r_t$ may have a negative drift when $K_t = 1$ but has a positive drift when $K_t \geq 2$ so we compare the fractions of time $K_t = 1$ and $K_t \geq 2$. We start by proving the result when

$$\phi_{aa} > \max\left(\frac{2\phi_{ba}}{5}, \phi_{ba} - \frac{\phi_{ab}}{6}\right) \quad \text{and} \quad \phi_{bb} = 0.$$

In this case, the drift when $K_t = 0$ is nonnegative (the worst case is drift = 0) so we shall think of state 0 as a "neutral" state, state 1 as a "bad" state and states in $I_2$ as "good" states. We will construct a Markov process $X_t$ with state space $\{0, 1, 2\}$ such that

(18) $$\frac{\limsup_t P(K_t = 1)}{\liminf_t P(K_t \geq 2)} \leq \frac{\lim_t P(X_t = 1)}{\lim_t P(X_t = 2)}.$$

Noticing that:

1. $K_t$ jumps from 1 to $I_2$ at rate at least $\phi_{ab} + \phi_{ba}$,
2. $K_t$ jumps from $I_2$ to $\{0, 1\}$ at rate at most $\phi_{aa} + \phi_{ab}/2 + \phi_{ba}$,

and having in mind that 0 is neutral, 1 is a bad state and 2 is a good state, we impose that:

1. $X_t$ cannot jump from 0 to 2 but jumps from 0 to 1 at some positive rate,

2. $X_t$ jumps from 1 to 2 at rate $\phi_{ab} + \phi_{ba}$,
3. $X_t$ jumps from 2 to $\{0,1\}$ at rate $\phi_{aa} + \phi_{ab}/2 + \phi_{ba}$.

The process $X_t$ describes a "worst case" scenario in which the time spent in bad states is rather large and the time spent in good states is rather small compared to $K_t$ so that (18) holds (note that the time spent in state 0 is unimportant in this inequality). Moreover, regardless of the specific values of the other transition rates of the process $X_t$, the conditions above imply that the stationary distribution $\pi = (\pi_0, \pi_1, \pi_2)$ of $X_t$ satisfies

$$0 \cdot \pi_0 + (\phi_{ab} + \phi_{ba}) \cdot \pi_1 - (\phi_{aa} + \phi_{ab}/2 + \phi_{ba}) \cdot \pi_2 = 0.$$

In particular, we can conclude that

$$\frac{\limsup_t P(K_t = 1)}{\liminf_t P(K_t \geq 2)} \leq \frac{\limsup_t P(X_t = 1)}{\liminf_t P(X_t = 2)}$$
$$= \frac{\lim_t P(X_t = 1)}{\lim_t P(X_t = 2)} = \frac{\pi_1}{\pi_2} = \frac{\phi_{aa} + \phi_{ab}/2 + \phi_{ba}}{\phi_{ab} + \phi_{ba}}$$

which we rewrite as

$$(\phi_{ab} + \phi_{ba}) \times \limsup_{t \to \infty} P(K_t = 1)$$
(19)
$$\leq (\phi_{aa} + \phi_{ab}/2 + \phi_{ba}) \times \liminf_{t \to \infty} P(K_t \geq 2).$$

To conclude the proof, we now distinguish two cases.

1. If $\phi_{aa} \leq (1/6)\phi_{ab}$ then (19) implies that

$$\limsup_{t \to \infty} P(K_t = 1) \leq \liminf_{t \to \infty} P(K_t \geq 2).$$

Moreover, the drift when $K_t \geq 2$ is larger than $4\phi_{aa} - \phi_{ba}$ and the rate at which the gap process $K_t$ jumps from 0 is strictly positive, so there exists $C_{20} > 0$ such that

$$\liminf_{t \to \infty} \frac{r_t}{t} \geq C_{20}(\phi_{aa} - \phi_{ba}) \limsup_{t \to \infty} P(K_t = 1)$$
$$+ C_{20}(4\phi_{aa} - \phi_{ba}) \liminf_{t \to \infty} P(K_t \geq 2)$$
$$\geq C_{20}(5\phi_{aa} - 2\phi_{ba}) \limsup_{t \to \infty} P(K_t = 1) > 0$$

since $\phi_{aa} > (2/5)\phi_{ba}$.

2. If $\phi_{aa} > (1/6)\phi_{ab}$, then, by using the fact that $\phi_{aa} \leq \phi_{ba}$ together with (19), we obtain

$$\limsup_{t \to \infty} P(K_t = 1) \leq 2 \times \liminf_{t \to \infty} P(K_t \geq 2).$$

Moreover, the drift when $K_t \geq 2$ is larger than $\phi_{aa} + (1/2)\phi_{ab} - \phi_{ba}$ and the rate at which the gap process $K_t$ jumps from 0 is strictly positive, so there exists $C_{21} > 0$ such that

$$\liminf_{t \to \infty} \frac{r_t}{t} \geq C_{21}(\phi_{aa} - \phi_{ba}) \limsup_{t \to \infty} P(K_t = 1)$$
$$+ C_{21}(\phi_{aa} + (1/2)\phi_{ab} - \phi_{ba}) \liminf_{t \to \infty} P(K_t \geq 2)$$
$$\geq C_{21}(3/2)(\phi_{aa} + (1/6)\phi_{ab} - \phi_{ba}) \limsup_{t \to \infty} P(K_t = 1) > 0$$

since $\phi_{aa} > \phi_{ba} - (1/6)\phi_{ab}$.

This establishes the result in the case when $\phi_{bb} = 0$. In addition, our proof implies that, properly rescaled in space and time, the process can be embedded in a supercritical oriented percolation process (we omit the details of the proof) which allows us to apply a perturbation argument similar to the one in the proof of Theorem 1 and extend the result to a set of parameters where $\phi_{bb}$ is small but positive. This completes the proof of Theorem 6.

**Acknowledgment.** The authors thank an anonymous referee for his/her numerous suggestions to improve some of our results, and three other referees for additional interesting comments.

School of Mathematical and
Statistical Sciences
Arizona State University
PO Box 871804
Tempe, Arizona 85287-1804
USA
E-mail: lanchier@math.asu.edu

Department of Ecology, Evolution
and Behavior
University of Minnesota
1987 Upper Buford Circle
St. Paul, Minnesota 55108
USA
E-mail: neuha001@umn.edu